\documentclass[final,onefignum,onetabnum]{siuro250211}

\usepackage{tikz}
\usepackage{pgfplots}
\pgfplotsset{compat=1.18}
\usepackage{amsfonts}
\usepackage{amsmath}
\usepackage{amssymb}
\usepackage{graphicx}
\usepackage{epstopdf}
\usepackage{float}
\usepackage{algorithmic}
\usepackage{algorithm}
\usepackage{caption}
\usepackage{subcaption}
\usepackage{booktabs}
\usepackage{bm}
\usetikzlibrary{plotmarks}
\ifpdf
  \DeclareGraphicsExtensions{.eps,.pdf,.png,.jpg}
\else
  \DeclareGraphicsExtensions{.eps}
\fi

\usepackage{enumitem}
\setlist[enumerate]{leftmargin=.5in}
\setlist[itemize]{leftmargin=.5in}

\newsiamremark{remark}{Remark}
\newsiamremark{hypothesis}{Hypothesis}
\crefname{hypothesis}{Hypothesis}{Hypotheses}
\newsiamthm{claim}{Claim}
\newsiamremark{fact}{Fact}
\crefname{fact}{Fact}{Facts}

\headers{Kronecker Approximations of Covariance Matrices for Solving Inverse Problems}{S. Sarkar}

\title{Kronecker Approximations of Covariance Matrices for Solving Inverse Problems}

\author{Srijon Sarkar\thanks{Department of Mathematics, Emory University
  (\email{srijon.sarkar@emory.edu}).}}

\dedication{\small\textit{Project advisor: Julianne Chung}}

\usepackage{amsopn}

\DeclareMathOperator{\vvec}{vec}
\DeclareMathOperator{\mat}{mat}

\ifpdf
\hypersetup{
  pdftitle={Kronecker Approximations of Covariance Matrices for Solving Inverse Problems},
  pdfauthor={S. Sarkar}
}
\fi

\begin{document}

\maketitle

\begin{abstract}
Inverse problems arise in many applications where the task is to recover unknown parameters from indirect or noisy observations. In contrast to forward problems, where outputs are directly computed from known inputs, inverse problems are designed to infer the underlying variables from the observed effect.  Inverse problems are often ill-posed, so regularization is needed to compute reasonable solutions. One form of regularization incorporates prior knowledge about the unknown parameters in the form of prior distributions.  However, in many settings, large prior covariance matrices pose computational bottlenecks for solving inverse problems.  To address this, we study approximations of large covariance matrices
via structured Kronecker--Toeplitz products and demonstrate how these approximations can be used to solve inverse problems.
We exploit the Block-Toeplitz-Toeplitz-Block structure induced by the covariance
kernel to construct compact approximations that are orders of magnitude cheaper
to store and invert. These approximations are deployed as prior covariance surrogates in Bayesian
inverse problems solved via generalized hybrid iterative methods.
We demonstrate that the Kronecker--Toeplitz approximation attains lower relative approximation error than existing matrix approximations, exhibits rapid error decay with rank, and enables efficient matrix--vector products.
Applications to 2D seismic tomography reconstruction problems illustrate the practical utility of these methods.
\end{abstract}

\section{Introduction}
\label{sec:intro}
Inverse problems arise in various scenarios across our everyday lives, from bioimaging to geology to aerospace, and the signals they involve affect major studies and environmental decisions. Given indirect observations and knowledge of the forward operator that produced them, the goal of an inverse problem is to recover the unknown parameters describing the underlying process or phenomenon. This is in contrast to forward problems, where these parameters are already known and the outputs are directly computable from the action of the measurement operator on the known inputs. Often the observations are polluted by noise or perturbations, e.g., from erroneous recordings or from nature. We are interested in inverse problems where, given a noisy real-world observation $\mathbf{b}$ and the forward operator $\textbf{A}$, we want to recover the true solution $\mathbf{x}$ in the linear system
\begin{equation}
  \mathbf{b} = \mathbf{A}\mathbf{x}_\text{true} + \bm{\epsilon},\label{eq:noisyinv}  
\end{equation}
where $\mathbf{x}_\text{true} \in \mathbb{R}^{n^2}$ is the unknown image on an $n \times n$ grid stored as a vector, i.e., $\mathbf{x} = \vvec(\mathbf{X})$ with $\mathbf{X} \in \mathbb{R}^{n \times n}$, $\mathbf{A} \in \mathbb{R}^{m\times n^2}$ is the forward operator, $\mathbf{b}\in\mathbb{R}^m$ contains the observations, and $\bm{\epsilon} \in \mathbb{R}^m$, $\bm{\epsilon} \sim \mathcal{N}(\mathbf{0},\mathbf{R})$, represents measurement error or noise. Here $\mathbf{R} \in \mathbb{R}^{m \times m}$ is a symmetric positive definite (SPD) noise covariance matrix, which is often assumed to be diagonal so that operations with $\mathbf{R}$ are inexpensive. The noise term $\bm{\epsilon}$ hinders accurate reconstruction of real-world outcomes. Further, inverse problems are typically \emph{ill-posed} by Hadamard's criteria~\cite{hadamard1902}. In other words, the singular values of $\mathbf{A}$ decay to and cluster at zero, resulting in noise amplification that dominates the reconstruction.


To prevent this noise propagation, regularization is used. There are different ways of applying regularization. In this paper, we consider variational regularization, specifically Tikhonov regularization. We assume a Gaussian prior for the solution, i.e., $\mathbf{x} \sim \mathcal{N}(0, \lambda^{-2}\mathbf{Q})$, where $\lambda > 0$ is the regularization parameter and $\mathbf{Q} \in \mathbb{R}^{n^2 \times n^2}$ is an SPD covariance matrix. Taking the Cholesky factorization $\mathbf{Q} = \mathbf{L}\mathbf{L}^T$, where $\mathbf{L} \in \mathbb{R}^{n^2 \times n^2}$ is lower triangular, we obtain a general-form Tikhonov problem, 
\begin{equation}
\min_\mathbf{x} \|\textbf{A}\mathbf{x} - \mathbf{b}\|_{\textbf{R}^{-1}}^2 + \lambda^2\|\mathbf{x}\|_{\textbf{Q}^{-1}}^2 = \min_\mathbf{x} \left\| \begin{bmatrix}
\mathbf{R}^{-\frac12} \mathbf{A} \\
{\lambda }\mathbf{L}^{-1}
\end{bmatrix}\mathbf{x} - \begin{bmatrix}
\mathbf{R}^{-\frac12} \mathbf{b} \\
\mathbf{0}
\end{bmatrix}\right\|_2^2, 
\label{eq:tikhonov}    
\end{equation}
where $\|\mathbf{x}\|_{\mathbf{M}}^2 = \mathbf{x}^T \mathbf{M}\mathbf{x}$ for SPD matrix $\mathbf{M}$.

Here, the regularization parameter $\lambda$ balances the goodness-of-fit (e.g., underfitting or overfitting the data) against the solution smoothness (e.g., controlling the noise that enters the solution). There also exist more general loss functions, additional constraints, and priors that could be incorporated~\cite{benning2018}. A primary difficulty lies in choosing $\lambda$ before obtaining a solution. One approach that does not require an explicit choice of the regularization parameter $\lambda$ is iterative regularization, where a projected problem is solved at each iteration and the stopping iteration plays the role of the regularization parameter~\cite{gazzola2019}. Hybrid projection methods pair both variational and iterative approaches to combine prior-informed, regularization techniques with reduced subspace frameworks for large-scale inverse problems that we visit later~\cite{chung2024}.

\begin{figure}[h!]
        \centering
        \resizebox{0.8\textwidth}{!}{
%
%
\begin{tikzpicture}

\begin{axis}[%
width=3.557in,
height=3.557in,
at={(2.167in,0.917in)},
scale only axis,
point meta min=-1308466115.4517,
point meta max=1218607942.29762,
axis on top,
xmin=0.5,
xmax=32.5,
tick align=outside,
y dir=reverse,
ymin=0.5,
ymax=32.5,
axis line style={draw=none},
ticks=none,
title style={font=\bfseries},
title={No Regularization}
]
\addplot [forget plot] graphics [xmin=0.5, xmax=32.5, ymin=0.5, ymax=32.5] {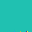};
\end{axis}

\begin{axis}[%
width=3.557in,
height=3.557in,
at={(6.847in,0.917in)},
scale only axis,
point meta min=-0.02012565090406,
point meta max=0.973905673724902,
axis on top,
xmin=0.5,
xmax=32.5,
tick align=outside,
y dir=reverse,
ymin=0.5,
ymax=32.5,
axis line style={draw=none},
ticks=none,
title style={font=\bfseries},
title={Regularized Solution}
]
\addplot [forget plot] graphics [xmin=0.5, xmax=32.5, ymin=0.5, ymax=32.5] {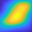};
\end{axis}

\begin{axis}[%
width=3.557in,
height=3.557in,
at={(11.527in,0.917in)},
scale only axis,
point meta min=0.0101607142558508,
point meta max=1,
axis on top,
xmin=0.5,
xmax=32.5,
tick align=outside,
y dir=reverse,
ymin=0.5,
ymax=32.5,
axis line style={draw=none},
ticks=none,
title style={font=\bfseries},
title={True Solution}
]
\addplot [forget plot] graphics [xmin=0.5, xmax=32.5, ymin=0.5, ymax=32.5] {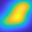};
\end{axis}
\end{tikzpicture}%
}
        \caption{Seismic tomography reconstruction of the slowness field on a $32\times32$ grid. Left: an unregularized least-squares solution of \eqref{eq:noisyinv}, where the small singular values of $\mathbf{A}$ amplify the noise leading to a poor reconstruction of the image. Center: a reconstruction with the Tikhonov regularization scheme \eqref{eq:tikhonov} using a Mat\'{e}rn prior covariance kernel. Right: the true slowness field. This example illustrates that Tikhonov regularization with the Mat\'{e}rn prior overcomes the noise amplification that occurs without regularization and can be used to recover an accurate, regularized solution.}
        \label{fig:example_noise}
\end{figure}

In \Cref{fig:example_noise}, we illustrate the importance of regularization for a seismic tomography problem.  We observe a big difference between the reconstruction of the solution when noise is not mitigated \cref{eq:noisyinv} compared to when it is \cref{eq:tikhonov}.

Incorporating well-informed prior knowledge, e.g., via the prior covariance matrix, can be crucial to the reconstruction quality of the solution.  However, $\mathbf{Q}$ is often large and dense, presenting one of the main computational bottlenecks of large-scale computational inverse problems. We address this challenge by proposing an efficient approximation of the prior covariance matrix to replace the original dense $\mathbf{Q}$. 

The paper is organized as follows: \Cref{sec:bg} covers the essential definitions and relevant formulas to preface our work. \Cref{sec:approx} introduces the Kronecker-Toeplitz approximation and compares it with existing approaches. \Cref{sec:inverse} formulates the Bayesian inverse problem and describes various iterative frameworks for solving \eqref{eq:tikhonov}, including Krylov subspace methods. \Cref{sec:formulations} presents different ways that our approximation can be deployed, and \Cref{sec:applications} provides numerical comparisons. Conclusions and future directions are provided in \Cref{sec:conclusions}.

\section{Background}
\label{sec:bg}

\begin{figure}[h!]
    \centering
    \begin{subfigure}[b]{0.31\linewidth}
        \resizebox{\linewidth}{!}{\input{pics/kernel}}
    \end{subfigure}\hfill
    \begin{subfigure}[b]{0.34\linewidth}
        \resizebox{\linewidth}{!}{
%
%
\begin{tikzpicture}

\begin{axis}[%
width=4.01in,
height=3.566in,
at={(0.673in,0.481in)},
scale only axis,
point meta min=0.7,
point meta max=1,
axis on top,
xmin=0.5,
xmax=1024.5,
y dir=reverse,
xtick={0, 200, 400, 600, 800, 1000},
ymin=0.5,
ymax=1024.5,
axis background/.style={fill=white},
colormap={mymap}{[1pt] rgb(0pt)=(0.2422,0.1504,0.6603); rgb(1pt)=(0.2444,0.1534,0.6728); rgb(2pt)=(0.2464,0.1569,0.6847); rgb(3pt)=(0.2484,0.1607,0.6961); rgb(4pt)=(0.2503,0.1648,0.7071); rgb(5pt)=(0.2522,0.1689,0.7179); rgb(6pt)=(0.254,0.1732,0.7286); rgb(7pt)=(0.2558,0.1773,0.7393); rgb(8pt)=(0.2576,0.1814,0.7501); rgb(9pt)=(0.2594,0.1854,0.761); rgb(11pt)=(0.2628,0.1932,0.7828); rgb(12pt)=(0.2645,0.1972,0.7937); rgb(13pt)=(0.2661,0.2011,0.8043); rgb(14pt)=(0.2676,0.2052,0.8148); rgb(15pt)=(0.2691,0.2094,0.8249); rgb(16pt)=(0.2704,0.2138,0.8346); rgb(17pt)=(0.2717,0.2184,0.8439); rgb(18pt)=(0.2729,0.2231,0.8528); rgb(19pt)=(0.274,0.228,0.8612); rgb(20pt)=(0.2749,0.233,0.8692); rgb(21pt)=(0.2758,0.2382,0.8767); rgb(22pt)=(0.2766,0.2435,0.884); rgb(23pt)=(0.2774,0.2489,0.8908); rgb(24pt)=(0.2781,0.2543,0.8973); rgb(25pt)=(0.2788,0.2598,0.9035); rgb(26pt)=(0.2794,0.2653,0.9094); rgb(27pt)=(0.2798,0.2708,0.915); rgb(28pt)=(0.2802,0.2764,0.9204); rgb(29pt)=(0.2806,0.2819,0.9255); rgb(30pt)=(0.2809,0.2875,0.9305); rgb(31pt)=(0.2811,0.293,0.9352); rgb(32pt)=(0.2813,0.2985,0.9397); rgb(33pt)=(0.2814,0.304,0.9441); rgb(34pt)=(0.2814,0.3095,0.9483); rgb(35pt)=(0.2813,0.315,0.9524); rgb(36pt)=(0.2811,0.3204,0.9563); rgb(37pt)=(0.2809,0.3259,0.96); rgb(38pt)=(0.2807,0.3313,0.9636); rgb(39pt)=(0.2803,0.3367,0.967); rgb(40pt)=(0.2798,0.3421,0.9702); rgb(41pt)=(0.2791,0.3475,0.9733); rgb(42pt)=(0.2784,0.3529,0.9763); rgb(43pt)=(0.2776,0.3583,0.9791); rgb(44pt)=(0.2766,0.3638,0.9817); rgb(45pt)=(0.2754,0.3693,0.984); rgb(46pt)=(0.2741,0.3748,0.9862); rgb(47pt)=(0.2726,0.3804,0.9881); rgb(48pt)=(0.271,0.386,0.9898); rgb(49pt)=(0.2691,0.3916,0.9912); rgb(50pt)=(0.267,0.3973,0.9924); rgb(51pt)=(0.2647,0.403,0.9935); rgb(52pt)=(0.2621,0.4088,0.9946); rgb(53pt)=(0.2591,0.4145,0.9955); rgb(54pt)=(0.2556,0.4203,0.9965); rgb(55pt)=(0.2517,0.4261,0.9974); rgb(56pt)=(0.2473,0.4319,0.9983); rgb(57pt)=(0.2424,0.4378,0.9991); rgb(58pt)=(0.2369,0.4437,0.9996); rgb(59pt)=(0.2311,0.4497,0.9995); rgb(60pt)=(0.225,0.4559,0.9985); rgb(61pt)=(0.2189,0.462,0.9968); rgb(62pt)=(0.2128,0.4682,0.9948); rgb(63pt)=(0.2066,0.4743,0.9926); rgb(64pt)=(0.2006,0.4803,0.9906); rgb(65pt)=(0.195,0.4861,0.9887); rgb(66pt)=(0.1903,0.4919,0.9867); rgb(67pt)=(0.1869,0.4975,0.9844); rgb(68pt)=(0.1847,0.503,0.9819); rgb(69pt)=(0.1831,0.5084,0.9793); rgb(70pt)=(0.1818,0.5138,0.9766); rgb(71pt)=(0.1806,0.5191,0.9738); rgb(72pt)=(0.1795,0.5244,0.9709); rgb(73pt)=(0.1785,0.5296,0.9677); rgb(74pt)=(0.1778,0.5349,0.9641); rgb(75pt)=(0.1773,0.5401,0.9602); rgb(76pt)=(0.1768,0.5452,0.956); rgb(77pt)=(0.1764,0.5504,0.9516); rgb(78pt)=(0.1755,0.5554,0.9473); rgb(79pt)=(0.174,0.5605,0.9432); rgb(80pt)=(0.1716,0.5655,0.9393); rgb(81pt)=(0.1686,0.5705,0.9357); rgb(82pt)=(0.1649,0.5755,0.9323); rgb(83pt)=(0.161,0.5805,0.9289); rgb(84pt)=(0.1573,0.5854,0.9254); rgb(85pt)=(0.154,0.5902,0.9218); rgb(86pt)=(0.1513,0.595,0.9182); rgb(87pt)=(0.1492,0.5997,0.9147); rgb(88pt)=(0.1475,0.6043,0.9113); rgb(89pt)=(0.1461,0.6089,0.908); rgb(90pt)=(0.1446,0.6135,0.905); rgb(91pt)=(0.1429,0.618,0.9022); rgb(92pt)=(0.1408,0.6226,0.8998); rgb(93pt)=(0.1383,0.6272,0.8975); rgb(94pt)=(0.1354,0.6317,0.8953); rgb(95pt)=(0.1321,0.6363,0.8932); rgb(96pt)=(0.1288,0.6408,0.891); rgb(97pt)=(0.1253,0.6453,0.8887); rgb(98pt)=(0.1219,0.6497,0.8862); rgb(99pt)=(0.1185,0.6541,0.8834); rgb(100pt)=(0.1152,0.6584,0.8804); rgb(101pt)=(0.1119,0.6627,0.877); rgb(102pt)=(0.1085,0.6669,0.8734); rgb(103pt)=(0.1048,0.671,0.8695); rgb(104pt)=(0.1009,0.675,0.8653); rgb(105pt)=(0.0964,0.6789,0.8609); rgb(106pt)=(0.0914,0.6828,0.8562); rgb(107pt)=(0.0855,0.6865,0.8513); rgb(108pt)=(0.0789,0.6902,0.8462); rgb(109pt)=(0.0713,0.6938,0.8409); rgb(110pt)=(0.0628,0.6972,0.8355); rgb(111pt)=(0.0535,0.7006,0.8299); rgb(112pt)=(0.0433,0.7039,0.8242); rgb(113pt)=(0.0328,0.7071,0.8183); rgb(114pt)=(0.0234,0.7103,0.8124); rgb(115pt)=(0.0155,0.7133,0.8064); rgb(116pt)=(0.0091,0.7163,0.8003); rgb(117pt)=(0.0046,0.7192,0.7941); rgb(118pt)=(0.0019,0.722,0.7878); rgb(119pt)=(0.0009,0.7248,0.7815); rgb(120pt)=(0.0018,0.7275,0.7752); rgb(121pt)=(0.0046,0.7301,0.7688); rgb(122pt)=(0.0094,0.7327,0.7623); rgb(123pt)=(0.0162,0.7352,0.7558); rgb(124pt)=(0.0253,0.7376,0.7492); rgb(125pt)=(0.0369,0.74,0.7426); rgb(126pt)=(0.0504,0.7423,0.7359); rgb(127pt)=(0.0638,0.7446,0.7292); rgb(128pt)=(0.077,0.7468,0.7224); rgb(129pt)=(0.0899,0.7489,0.7156); rgb(130pt)=(0.1023,0.751,0.7088); rgb(131pt)=(0.1141,0.7531,0.7019); rgb(132pt)=(0.1252,0.7552,0.695); rgb(133pt)=(0.1354,0.7572,0.6881); rgb(134pt)=(0.1448,0.7593,0.6812); rgb(135pt)=(0.1532,0.7614,0.6741); rgb(136pt)=(0.1609,0.7635,0.6671); rgb(137pt)=(0.1678,0.7656,0.6599); rgb(138pt)=(0.1741,0.7678,0.6527); rgb(139pt)=(0.1799,0.7699,0.6454); rgb(140pt)=(0.1853,0.7721,0.6379); rgb(141pt)=(0.1905,0.7743,0.6303); rgb(142pt)=(0.1954,0.7765,0.6225); rgb(143pt)=(0.2003,0.7787,0.6146); rgb(144pt)=(0.2061,0.7808,0.6065); rgb(145pt)=(0.2118,0.7828,0.5983); rgb(146pt)=(0.2178,0.7849,0.5899); rgb(147pt)=(0.2244,0.7869,0.5813); rgb(148pt)=(0.2318,0.7887,0.5725); rgb(149pt)=(0.2401,0.7905,0.5636); rgb(150pt)=(0.2491,0.7922,0.5546); rgb(151pt)=(0.2589,0.7937,0.5454); rgb(152pt)=(0.2695,0.7951,0.536); rgb(153pt)=(0.2809,0.7964,0.5266); rgb(154pt)=(0.2929,0.7975,0.517); rgb(155pt)=(0.3052,0.7985,0.5074); rgb(156pt)=(0.3176,0.7994,0.4975); rgb(157pt)=(0.3301,0.8002,0.4876); rgb(158pt)=(0.3424,0.8009,0.4774); rgb(159pt)=(0.3548,0.8016,0.4669); rgb(160pt)=(0.3671,0.8021,0.4563); rgb(161pt)=(0.3795,0.8026,0.4454); rgb(162pt)=(0.3921,0.8029,0.4344); rgb(163pt)=(0.405,0.8031,0.4233); rgb(164pt)=(0.4184,0.803,0.4122); rgb(165pt)=(0.4322,0.8028,0.4013); rgb(166pt)=(0.4463,0.8024,0.3904); rgb(167pt)=(0.4608,0.8018,0.3797); rgb(168pt)=(0.4753,0.8011,0.3691); rgb(169pt)=(0.4899,0.8002,0.3586); rgb(170pt)=(0.5044,0.7993,0.348); rgb(171pt)=(0.5187,0.7982,0.3374); rgb(172pt)=(0.5329,0.797,0.3267); rgb(173pt)=(0.547,0.7957,0.3159); rgb(175pt)=(0.5748,0.7929,0.2941); rgb(176pt)=(0.5886,0.7913,0.2833); rgb(177pt)=(0.6024,0.7896,0.2726); rgb(178pt)=(0.6161,0.7878,0.2622); rgb(179pt)=(0.6297,0.7859,0.2521); rgb(180pt)=(0.6433,0.7839,0.2423); rgb(181pt)=(0.6567,0.7818,0.2329); rgb(182pt)=(0.6701,0.7796,0.2239); rgb(183pt)=(0.6833,0.7773,0.2155); rgb(184pt)=(0.6963,0.775,0.2075); rgb(185pt)=(0.7091,0.7727,0.1998); rgb(186pt)=(0.7218,0.7703,0.1924); rgb(187pt)=(0.7344,0.7679,0.1852); rgb(188pt)=(0.7468,0.7654,0.1782); rgb(189pt)=(0.759,0.7629,0.1717); rgb(190pt)=(0.771,0.7604,0.1658); rgb(191pt)=(0.7829,0.7579,0.1608); rgb(192pt)=(0.7945,0.7554,0.157); rgb(193pt)=(0.806,0.7529,0.1546); rgb(194pt)=(0.8172,0.7505,0.1535); rgb(195pt)=(0.8281,0.7481,0.1536); rgb(196pt)=(0.8389,0.7457,0.1546); rgb(197pt)=(0.8495,0.7435,0.1564); rgb(198pt)=(0.86,0.7413,0.1587); rgb(199pt)=(0.8703,0.7392,0.1615); rgb(200pt)=(0.8804,0.7372,0.165); rgb(201pt)=(0.8903,0.7353,0.1695); rgb(202pt)=(0.9,0.7336,0.1749); rgb(203pt)=(0.9093,0.7321,0.1815); rgb(204pt)=(0.9184,0.7308,0.189); rgb(205pt)=(0.9272,0.7298,0.1973); rgb(206pt)=(0.9357,0.729,0.2061); rgb(207pt)=(0.944,0.7285,0.2151); rgb(208pt)=(0.9523,0.7284,0.2237); rgb(209pt)=(0.9606,0.7285,0.2312); rgb(210pt)=(0.9689,0.7292,0.2373); rgb(211pt)=(0.977,0.7304,0.2418); rgb(212pt)=(0.9842,0.733,0.2446); rgb(213pt)=(0.99,0.7365,0.2429); rgb(214pt)=(0.9946,0.7407,0.2394); rgb(215pt)=(0.9966,0.7458,0.2351); rgb(216pt)=(0.9971,0.7513,0.2309); rgb(217pt)=(0.9972,0.7569,0.2267); rgb(218pt)=(0.9971,0.7626,0.2224); rgb(219pt)=(0.9969,0.7683,0.2181); rgb(220pt)=(0.9966,0.774,0.2138); rgb(221pt)=(0.9962,0.7798,0.2095); rgb(222pt)=(0.9957,0.7856,0.2053); rgb(223pt)=(0.9949,0.7915,0.2012); rgb(224pt)=(0.9938,0.7974,0.1974); rgb(225pt)=(0.9923,0.8034,0.1939); rgb(226pt)=(0.9906,0.8095,0.1906); rgb(227pt)=(0.9885,0.8156,0.1875); rgb(228pt)=(0.9861,0.8218,0.1846); rgb(229pt)=(0.9835,0.828,0.1817); rgb(230pt)=(0.9807,0.8342,0.1787); rgb(231pt)=(0.9778,0.8404,0.1757); rgb(232pt)=(0.9748,0.8467,0.1726); rgb(233pt)=(0.972,0.8529,0.1695); rgb(234pt)=(0.9694,0.8591,0.1665); rgb(235pt)=(0.9671,0.8654,0.1636); rgb(236pt)=(0.9651,0.8716,0.1608); rgb(237pt)=(0.9634,0.8778,0.1582); rgb(238pt)=(0.9619,0.884,0.1557); rgb(239pt)=(0.9608,0.8902,0.1532); rgb(240pt)=(0.9601,0.8963,0.1507); rgb(241pt)=(0.9596,0.9023,0.148); rgb(242pt)=(0.9595,0.9084,0.145); rgb(243pt)=(0.9597,0.9143,0.1418); rgb(244pt)=(0.9601,0.9203,0.1382); rgb(245pt)=(0.9608,0.9262,0.1344); rgb(246pt)=(0.9618,0.932,0.1304); rgb(247pt)=(0.9629,0.9379,0.1261); rgb(248pt)=(0.9642,0.9437,0.1216); rgb(249pt)=(0.9657,0.9494,0.1168); rgb(250pt)=(0.9674,0.9552,0.1116); rgb(251pt)=(0.9692,0.9609,0.1061); rgb(252pt)=(0.9711,0.9667,0.1001); rgb(253pt)=(0.973,0.9724,0.0938); rgb(254pt)=(0.9749,0.9782,0.0872); rgb(255pt)=(0.9769,0.9839,0.0805)},
colorbar
]
\addplot [forget plot] graphics [xmin=0.5, xmax=1024.5, ymin=0.5, ymax=1024.5] {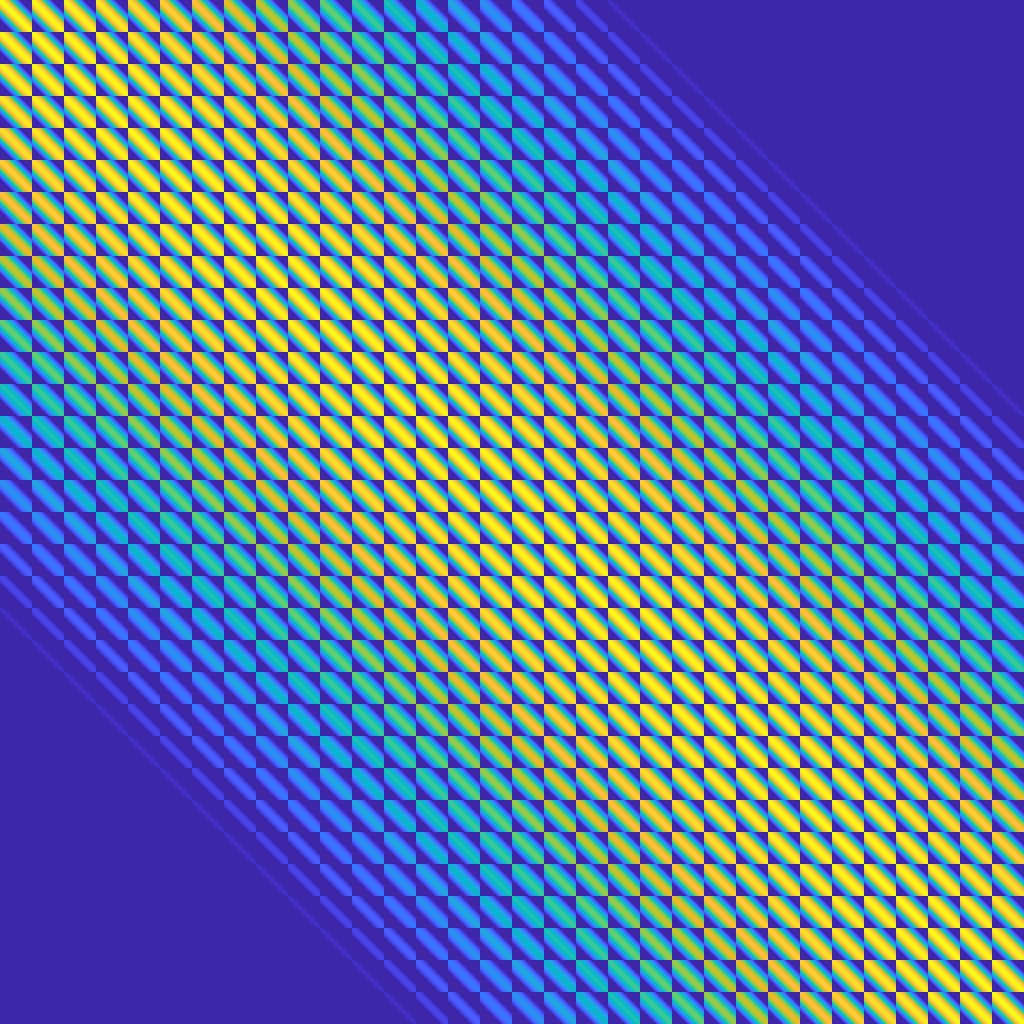};
\end{axis}
\end{tikzpicture}%
}
    \end{subfigure}\hfill
    \begin{subfigure}[b]{0.31\linewidth}
        \resizebox{\linewidth}{!}{
%
%
\begin{tikzpicture}

\begin{axis}[%
width=4.047in,
height=3.566in,
at={(0.679in,0.481in)},
scale only axis,
point meta min=0.7,
point meta max=1,
axis on top,
xmin=0.5,
xmax=32.5,
y dir=reverse,
ymin=0.5,
ymax=32.5,
axis background/.style={fill=white},
title style={font=\bfseries},
title={PSF (first column)},
colormap={mymap}{[1pt] rgb(0pt)=(0.2422,0.1504,0.6603); rgb(1pt)=(0.2444,0.1534,0.6728); rgb(2pt)=(0.2464,0.1569,0.6847); rgb(3pt)=(0.2484,0.1607,0.6961); rgb(4pt)=(0.2503,0.1648,0.7071); rgb(5pt)=(0.2522,0.1689,0.7179); rgb(6pt)=(0.254,0.1732,0.7286); rgb(7pt)=(0.2558,0.1773,0.7393); rgb(8pt)=(0.2576,0.1814,0.7501); rgb(9pt)=(0.2594,0.1854,0.761); rgb(11pt)=(0.2628,0.1932,0.7828); rgb(12pt)=(0.2645,0.1972,0.7937); rgb(13pt)=(0.2661,0.2011,0.8043); rgb(14pt)=(0.2676,0.2052,0.8148); rgb(15pt)=(0.2691,0.2094,0.8249); rgb(16pt)=(0.2704,0.2138,0.8346); rgb(17pt)=(0.2717,0.2184,0.8439); rgb(18pt)=(0.2729,0.2231,0.8528); rgb(19pt)=(0.274,0.228,0.8612); rgb(20pt)=(0.2749,0.233,0.8692); rgb(21pt)=(0.2758,0.2382,0.8767); rgb(22pt)=(0.2766,0.2435,0.884); rgb(23pt)=(0.2774,0.2489,0.8908); rgb(24pt)=(0.2781,0.2543,0.8973); rgb(25pt)=(0.2788,0.2598,0.9035); rgb(26pt)=(0.2794,0.2653,0.9094); rgb(27pt)=(0.2798,0.2708,0.915); rgb(28pt)=(0.2802,0.2764,0.9204); rgb(29pt)=(0.2806,0.2819,0.9255); rgb(30pt)=(0.2809,0.2875,0.9305); rgb(31pt)=(0.2811,0.293,0.9352); rgb(32pt)=(0.2813,0.2985,0.9397); rgb(33pt)=(0.2814,0.304,0.9441); rgb(34pt)=(0.2814,0.3095,0.9483); rgb(35pt)=(0.2813,0.315,0.9524); rgb(36pt)=(0.2811,0.3204,0.9563); rgb(37pt)=(0.2809,0.3259,0.96); rgb(38pt)=(0.2807,0.3313,0.9636); rgb(39pt)=(0.2803,0.3367,0.967); rgb(40pt)=(0.2798,0.3421,0.9702); rgb(41pt)=(0.2791,0.3475,0.9733); rgb(42pt)=(0.2784,0.3529,0.9763); rgb(43pt)=(0.2776,0.3583,0.9791); rgb(44pt)=(0.2766,0.3638,0.9817); rgb(45pt)=(0.2754,0.3693,0.984); rgb(46pt)=(0.2741,0.3748,0.9862); rgb(47pt)=(0.2726,0.3804,0.9881); rgb(48pt)=(0.271,0.386,0.9898); rgb(49pt)=(0.2691,0.3916,0.9912); rgb(50pt)=(0.267,0.3973,0.9924); rgb(51pt)=(0.2647,0.403,0.9935); rgb(52pt)=(0.2621,0.4088,0.9946); rgb(53pt)=(0.2591,0.4145,0.9955); rgb(54pt)=(0.2556,0.4203,0.9965); rgb(55pt)=(0.2517,0.4261,0.9974); rgb(56pt)=(0.2473,0.4319,0.9983); rgb(57pt)=(0.2424,0.4378,0.9991); rgb(58pt)=(0.2369,0.4437,0.9996); rgb(59pt)=(0.2311,0.4497,0.9995); rgb(60pt)=(0.225,0.4559,0.9985); rgb(61pt)=(0.2189,0.462,0.9968); rgb(62pt)=(0.2128,0.4682,0.9948); rgb(63pt)=(0.2066,0.4743,0.9926); rgb(64pt)=(0.2006,0.4803,0.9906); rgb(65pt)=(0.195,0.4861,0.9887); rgb(66pt)=(0.1903,0.4919,0.9867); rgb(67pt)=(0.1869,0.4975,0.9844); rgb(68pt)=(0.1847,0.503,0.9819); rgb(69pt)=(0.1831,0.5084,0.9793); rgb(70pt)=(0.1818,0.5138,0.9766); rgb(71pt)=(0.1806,0.5191,0.9738); rgb(72pt)=(0.1795,0.5244,0.9709); rgb(73pt)=(0.1785,0.5296,0.9677); rgb(74pt)=(0.1778,0.5349,0.9641); rgb(75pt)=(0.1773,0.5401,0.9602); rgb(76pt)=(0.1768,0.5452,0.956); rgb(77pt)=(0.1764,0.5504,0.9516); rgb(78pt)=(0.1755,0.5554,0.9473); rgb(79pt)=(0.174,0.5605,0.9432); rgb(80pt)=(0.1716,0.5655,0.9393); rgb(81pt)=(0.1686,0.5705,0.9357); rgb(82pt)=(0.1649,0.5755,0.9323); rgb(83pt)=(0.161,0.5805,0.9289); rgb(84pt)=(0.1573,0.5854,0.9254); rgb(85pt)=(0.154,0.5902,0.9218); rgb(86pt)=(0.1513,0.595,0.9182); rgb(87pt)=(0.1492,0.5997,0.9147); rgb(88pt)=(0.1475,0.6043,0.9113); rgb(89pt)=(0.1461,0.6089,0.908); rgb(90pt)=(0.1446,0.6135,0.905); rgb(91pt)=(0.1429,0.618,0.9022); rgb(92pt)=(0.1408,0.6226,0.8998); rgb(93pt)=(0.1383,0.6272,0.8975); rgb(94pt)=(0.1354,0.6317,0.8953); rgb(95pt)=(0.1321,0.6363,0.8932); rgb(96pt)=(0.1288,0.6408,0.891); rgb(97pt)=(0.1253,0.6453,0.8887); rgb(98pt)=(0.1219,0.6497,0.8862); rgb(99pt)=(0.1185,0.6541,0.8834); rgb(100pt)=(0.1152,0.6584,0.8804); rgb(101pt)=(0.1119,0.6627,0.877); rgb(102pt)=(0.1085,0.6669,0.8734); rgb(103pt)=(0.1048,0.671,0.8695); rgb(104pt)=(0.1009,0.675,0.8653); rgb(105pt)=(0.0964,0.6789,0.8609); rgb(106pt)=(0.0914,0.6828,0.8562); rgb(107pt)=(0.0855,0.6865,0.8513); rgb(108pt)=(0.0789,0.6902,0.8462); rgb(109pt)=(0.0713,0.6938,0.8409); rgb(110pt)=(0.0628,0.6972,0.8355); rgb(111pt)=(0.0535,0.7006,0.8299); rgb(112pt)=(0.0433,0.7039,0.8242); rgb(113pt)=(0.0328,0.7071,0.8183); rgb(114pt)=(0.0234,0.7103,0.8124); rgb(115pt)=(0.0155,0.7133,0.8064); rgb(116pt)=(0.0091,0.7163,0.8003); rgb(117pt)=(0.0046,0.7192,0.7941); rgb(118pt)=(0.0019,0.722,0.7878); rgb(119pt)=(0.0009,0.7248,0.7815); rgb(120pt)=(0.0018,0.7275,0.7752); rgb(121pt)=(0.0046,0.7301,0.7688); rgb(122pt)=(0.0094,0.7327,0.7623); rgb(123pt)=(0.0162,0.7352,0.7558); rgb(124pt)=(0.0253,0.7376,0.7492); rgb(125pt)=(0.0369,0.74,0.7426); rgb(126pt)=(0.0504,0.7423,0.7359); rgb(127pt)=(0.0638,0.7446,0.7292); rgb(128pt)=(0.077,0.7468,0.7224); rgb(129pt)=(0.0899,0.7489,0.7156); rgb(130pt)=(0.1023,0.751,0.7088); rgb(131pt)=(0.1141,0.7531,0.7019); rgb(132pt)=(0.1252,0.7552,0.695); rgb(133pt)=(0.1354,0.7572,0.6881); rgb(134pt)=(0.1448,0.7593,0.6812); rgb(135pt)=(0.1532,0.7614,0.6741); rgb(136pt)=(0.1609,0.7635,0.6671); rgb(137pt)=(0.1678,0.7656,0.6599); rgb(138pt)=(0.1741,0.7678,0.6527); rgb(139pt)=(0.1799,0.7699,0.6454); rgb(140pt)=(0.1853,0.7721,0.6379); rgb(141pt)=(0.1905,0.7743,0.6303); rgb(142pt)=(0.1954,0.7765,0.6225); rgb(143pt)=(0.2003,0.7787,0.6146); rgb(144pt)=(0.2061,0.7808,0.6065); rgb(145pt)=(0.2118,0.7828,0.5983); rgb(146pt)=(0.2178,0.7849,0.5899); rgb(147pt)=(0.2244,0.7869,0.5813); rgb(148pt)=(0.2318,0.7887,0.5725); rgb(149pt)=(0.2401,0.7905,0.5636); rgb(150pt)=(0.2491,0.7922,0.5546); rgb(151pt)=(0.2589,0.7937,0.5454); rgb(152pt)=(0.2695,0.7951,0.536); rgb(153pt)=(0.2809,0.7964,0.5266); rgb(154pt)=(0.2929,0.7975,0.517); rgb(155pt)=(0.3052,0.7985,0.5074); rgb(156pt)=(0.3176,0.7994,0.4975); rgb(157pt)=(0.3301,0.8002,0.4876); rgb(158pt)=(0.3424,0.8009,0.4774); rgb(159pt)=(0.3548,0.8016,0.4669); rgb(160pt)=(0.3671,0.8021,0.4563); rgb(161pt)=(0.3795,0.8026,0.4454); rgb(162pt)=(0.3921,0.8029,0.4344); rgb(163pt)=(0.405,0.8031,0.4233); rgb(164pt)=(0.4184,0.803,0.4122); rgb(165pt)=(0.4322,0.8028,0.4013); rgb(166pt)=(0.4463,0.8024,0.3904); rgb(167pt)=(0.4608,0.8018,0.3797); rgb(168pt)=(0.4753,0.8011,0.3691); rgb(169pt)=(0.4899,0.8002,0.3586); rgb(170pt)=(0.5044,0.7993,0.348); rgb(171pt)=(0.5187,0.7982,0.3374); rgb(172pt)=(0.5329,0.797,0.3267); rgb(173pt)=(0.547,0.7957,0.3159); rgb(175pt)=(0.5748,0.7929,0.2941); rgb(176pt)=(0.5886,0.7913,0.2833); rgb(177pt)=(0.6024,0.7896,0.2726); rgb(178pt)=(0.6161,0.7878,0.2622); rgb(179pt)=(0.6297,0.7859,0.2521); rgb(180pt)=(0.6433,0.7839,0.2423); rgb(181pt)=(0.6567,0.7818,0.2329); rgb(182pt)=(0.6701,0.7796,0.2239); rgb(183pt)=(0.6833,0.7773,0.2155); rgb(184pt)=(0.6963,0.775,0.2075); rgb(185pt)=(0.7091,0.7727,0.1998); rgb(186pt)=(0.7218,0.7703,0.1924); rgb(187pt)=(0.7344,0.7679,0.1852); rgb(188pt)=(0.7468,0.7654,0.1782); rgb(189pt)=(0.759,0.7629,0.1717); rgb(190pt)=(0.771,0.7604,0.1658); rgb(191pt)=(0.7829,0.7579,0.1608); rgb(192pt)=(0.7945,0.7554,0.157); rgb(193pt)=(0.806,0.7529,0.1546); rgb(194pt)=(0.8172,0.7505,0.1535); rgb(195pt)=(0.8281,0.7481,0.1536); rgb(196pt)=(0.8389,0.7457,0.1546); rgb(197pt)=(0.8495,0.7435,0.1564); rgb(198pt)=(0.86,0.7413,0.1587); rgb(199pt)=(0.8703,0.7392,0.1615); rgb(200pt)=(0.8804,0.7372,0.165); rgb(201pt)=(0.8903,0.7353,0.1695); rgb(202pt)=(0.9,0.7336,0.1749); rgb(203pt)=(0.9093,0.7321,0.1815); rgb(204pt)=(0.9184,0.7308,0.189); rgb(205pt)=(0.9272,0.7298,0.1973); rgb(206pt)=(0.9357,0.729,0.2061); rgb(207pt)=(0.944,0.7285,0.2151); rgb(208pt)=(0.9523,0.7284,0.2237); rgb(209pt)=(0.9606,0.7285,0.2312); rgb(210pt)=(0.9689,0.7292,0.2373); rgb(211pt)=(0.977,0.7304,0.2418); rgb(212pt)=(0.9842,0.733,0.2446); rgb(213pt)=(0.99,0.7365,0.2429); rgb(214pt)=(0.9946,0.7407,0.2394); rgb(215pt)=(0.9966,0.7458,0.2351); rgb(216pt)=(0.9971,0.7513,0.2309); rgb(217pt)=(0.9972,0.7569,0.2267); rgb(218pt)=(0.9971,0.7626,0.2224); rgb(219pt)=(0.9969,0.7683,0.2181); rgb(220pt)=(0.9966,0.774,0.2138); rgb(221pt)=(0.9962,0.7798,0.2095); rgb(222pt)=(0.9957,0.7856,0.2053); rgb(223pt)=(0.9949,0.7915,0.2012); rgb(224pt)=(0.9938,0.7974,0.1974); rgb(225pt)=(0.9923,0.8034,0.1939); rgb(226pt)=(0.9906,0.8095,0.1906); rgb(227pt)=(0.9885,0.8156,0.1875); rgb(228pt)=(0.9861,0.8218,0.1846); rgb(229pt)=(0.9835,0.828,0.1817); rgb(230pt)=(0.9807,0.8342,0.1787); rgb(231pt)=(0.9778,0.8404,0.1757); rgb(232pt)=(0.9748,0.8467,0.1726); rgb(233pt)=(0.972,0.8529,0.1695); rgb(234pt)=(0.9694,0.8591,0.1665); rgb(235pt)=(0.9671,0.8654,0.1636); rgb(236pt)=(0.9651,0.8716,0.1608); rgb(237pt)=(0.9634,0.8778,0.1582); rgb(238pt)=(0.9619,0.884,0.1557); rgb(239pt)=(0.9608,0.8902,0.1532); rgb(240pt)=(0.9601,0.8963,0.1507); rgb(241pt)=(0.9596,0.9023,0.148); rgb(242pt)=(0.9595,0.9084,0.145); rgb(243pt)=(0.9597,0.9143,0.1418); rgb(244pt)=(0.9601,0.9203,0.1382); rgb(245pt)=(0.9608,0.9262,0.1344); rgb(246pt)=(0.9618,0.932,0.1304); rgb(247pt)=(0.9629,0.9379,0.1261); rgb(248pt)=(0.9642,0.9437,0.1216); rgb(249pt)=(0.9657,0.9494,0.1168); rgb(250pt)=(0.9674,0.9552,0.1116); rgb(251pt)=(0.9692,0.9609,0.1061); rgb(252pt)=(0.9711,0.9667,0.1001); rgb(253pt)=(0.973,0.9724,0.0938); rgb(254pt)=(0.9749,0.9782,0.0872); rgb(255pt)=(0.9769,0.9839,0.0805)},
colorbar
]
\addplot [forget plot] graphics [xmin=0.5, xmax=32.5, ymin=0.5, ymax=32.5] {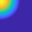};
\end{axis}
\end{tikzpicture}%
}
    \end{subfigure}
    \caption{We showcase 1D Mat\'ern kernels for different smoothness parameters and, further, a resulting covariance matrix for a 2D Mat\'ern and its corresponding PSF.}
    \label{fig:kernelQpsf}
\end{figure}

The covariance matrix $\textbf{Q}$ in \cref{eq:tikhonov} encodes our prior beliefs about $\mathbf{x}$ before observing any data. There are various approaches to select the prior covariance matrix $\mathbf{Q}$. Note that the general-form Tikhonov problem uses the inverse of the covariance matrix, often referred to as the precision matrix.
In many cases, computing the precision matrix is infeasible. This may occur, for example, when $\mathbf{Q}$ is very large and dense. However, efficient Matrix--Vector Products (MVPs) and operations may be possible. In this work we focus on 2D Mat\'{e}rn kernels, $k_\nu(r)$, defined as
\begin{equation}
\label{eq:matern}
  k_\nu(r) = \frac{2^{1-\nu}}{\Gamma(\nu)}
    \!\left(\frac{\sqrt{2\nu}\, r}{\ell}\right)^{\!\nu}
    K_\nu\!\left(\frac{\sqrt{2\nu}\, r}{\ell}\right),
\end{equation}
where $r$ denotes radial distance, $\nu$ the smoothness parameter, $\ell$ the correlation length, and $K_\nu(\cdot)$ the modified Bessel function of the second kind. As $\nu \to \infty$, this converges to the Gaussian (squared-exponential) kernel. Commonly used values are $\nu = \tfrac{1}{2}, \tfrac{3}{2}, \tfrac{5}{2}$.  In the left plot of \Cref{fig:kernelQpsf}, we provide 1D Mat\'{e}rn kernels for various choices of $\nu$, demonstrating the different decay behaviors. The parameter $\ell$ governs the rate of decay, where smaller $\ell$ gives correlations that decay over shorter distances. 

Let $\{\mathbf{t}_i\}_{i=1}^{n^2}$ be the nodes of a uniform $n \times n$ grid in $[0,1]^2$, ordered columnwise to match $\mathbf{x} = \vvec(\mathbf{X})$. Throughout, for a matrix $\mathbf{M}$, we denote its $j$th column by $\mathbf{m}_j$ and its $ij$th element by $m_{ij}$ (or $[\mathbf{M}]_{ij}$ when $\mathbf{M}$ carries a subscript). The prior covariance matrix $\mathbf{Q} \in \mathbb{R}^{n^2 \times n^2}$ is computed by applying the Mat\'{e}rn covariance function \eqref{eq:matern} to distances corresponding to every pair of grid nodes, so that the $ij$th element of the matrix is 
\begin{equation}
  q_{ij} = k_\nu\!\left(\|\mathbf{t}_i - \mathbf{t}_j\|_2\right).
\end{equation}
The matrix elements thus depend only on the Euclidean distance between any two grid nodes. In particular, $q_{ii} = k_\nu(0) = 1$ yields a unit diagonal (unit marginal variance), and $q_{ij} = q_{ji}$ since distance is symmetric. The smaller $\nu$ is, the rougher the field and hence the steeper the decay.  The larger $\ell$ is, the more correlated distant points are. Since the kernel is positive definite over $\mathbb{R}^2$, the obtained matrix $\mathbf{Q}$ is symmetric positive definite. In addition, because $k_\nu$ is a function of the difference between $\mathbf{t}_i$ and $\mathbf{t}_j$ only and the grid is uniform, the entries are constant along each block-diagonal, i.e., the matrix has \emph{Block-Toeplitz-Toeplitz-Block} (BTTB) form.


Mat\`{e}rn covariance matrices are highly structured, in that for regular grids, the matrix is BTTB. A square matrix is \emph{Toeplitz} if each descending or ascending diagonal, including the main diagonal, has the same value (constant) in all entries of that subdiagonal. Thus, a Toeplitz matrix can be defined by a single vector. A BTTB matrix is a 2D extension, where the matrix is partitioned into blocks, and both the arrangement of blocks and the individual blocks themselves follow this constant-diagonal structure. In our case, we will be dealing with an $n \times n$ block Toeplitz matrix with $n \times n$ Toeplitz blocks. The Mat\'{e}rn covariance matrix $\mathbf{Q}$ has BTTB structure since the kernel entries depend only on displacements, not absolute positions of the points in the domain. See the center plot of \Cref{fig:kernelQpsf} for a visualization of the BTTB structure. Recall that the \emph{Kronecker product} of matrices $\mathbf{A} \in \mathbb{R}^{n\times n}$ and $\mathbf{B} \in \mathbb{R}^{n\times n}$ is the $n^2 \times n^2$ block matrix \begin{equation}
  \mathbf{A} \otimes \mathbf{B} =
  \begin{bmatrix}
    a_{11}\mathbf{B} & \cdots & a_{1n}\mathbf{B} \\
    \vdots & \ddots & \vdots \\
    a_{n1}\mathbf{B} & \cdots & a_{nn}\mathbf{B}
  \end{bmatrix},
\end{equation}
whose $ij$-th block is the $n \times n$ matrix $a_{ij}\mathbf{B}$. Following~\cite{kamm2000}, $\vvec$ is an operator that maps a matrix to a vector by stacking its columns, i.e., for matrix $\!\mathbf{X} = [\mathbf{x}_1, \ldots, \mathbf{x}_n] \in \mathbb{R}^{n \times n}$,
\begin{equation}
\label{eq:vec}
  \vvec(\mathbf{X}) = \begin{bmatrix} \mathbf{x}_1 \\ \vdots \\ \mathbf{x}_n \end{bmatrix} \in \mathbb{R}^{n^2},
\end{equation}
and $\mat$ denotes the inverse operation.
By exploiting the following relationship for MVPs with Kronecker product matrices~\cite{kamm2000}, we have
\begin{equation}
\label{eq:kronvec}
  (\mathbf{A} \otimes \mathbf{B})\,\vvec(\mathbf{X}) = \vvec\!\left(\mathbf{B}\mathbf{X}\mathbf{A}^T\right).
\end{equation}
Thus, we can obtain an MVP with an $n^2 \times n^2$ Kronecker product matrix using only two matrix-matrix products among $n \times n$ matrices. In the problem of image deblurring, where BTTB matrices describe spatially invariant blurring, a matrix of such type is defined by its \emph{point spread function} (PSF), i.e., by the $n \times n$ image of the point source, representing how each pixel is smeared by its neighbors. We denote the PSF matrix by $\textbf{P}.$

Kronecker product approximations have been considered in many contexts~\cite{nagy1996, kammnagy1998, kilmer2022}, and similar approximations can be used for kernel matrices more generally. Kamm and Nagy~\cite{kamm2000} showed that for a banded BTTB matrix such as $\mathbf{Q}$, the factors $\mathbf{A}_k$ and $\mathbf{B}_k$ that provide the best Frobenius-norm approximation $\sum_{k=1}^{s}\mathbf{A}_k\otimes\mathbf{B}_k$ are banded Toeplitz matrices. These factors can be computed from the SVD of a \emph{weighted} PSF. In their paper, the Kronecker terms are further used to build an approximate SVD of $\mathbf{Q}$ whose small singular values, below a threshold chosen by generalized cross-validation, are set to one. The result is a preconditioner to accelerate iterative solvers, rather than an approximation of $\mathbf{Q}$ itself. 
In this work, we will consider Kronecker product approximations of large Mat\'{e}rn covariance matrices and use them within generalized hybrid iterative methods for solving inverse problems.

\section{The Kronecker--Toeplitz Approximation}
\label{sec:approx}

\begin{figure}[h!]
    \centering
    \begin{subfigure}[b]{0.33\linewidth}
        \resizebox{\linewidth}{!}{
%
%
\begin{tikzpicture}

\begin{axis}[%
width=4.01in,
height=3.566in,
at={(0.673in,0.481in)},
scale only axis,
point meta min=0.7,
point meta max=1,
axis on top,
xmin=0.5,
xmax=1024.5,
y dir=reverse,
xtick={0, 200, 400, 600, 800, 1000},
ymin=0.5,
ymax=1024.5,
axis background/.style={fill=white},
colormap={mymap}{[1pt] rgb(0pt)=(0.2422,0.1504,0.6603); rgb(1pt)=(0.2444,0.1534,0.6728); rgb(2pt)=(0.2464,0.1569,0.6847); rgb(3pt)=(0.2484,0.1607,0.6961); rgb(4pt)=(0.2503,0.1648,0.7071); rgb(5pt)=(0.2522,0.1689,0.7179); rgb(6pt)=(0.254,0.1732,0.7286); rgb(7pt)=(0.2558,0.1773,0.7393); rgb(8pt)=(0.2576,0.1814,0.7501); rgb(9pt)=(0.2594,0.1854,0.761); rgb(11pt)=(0.2628,0.1932,0.7828); rgb(12pt)=(0.2645,0.1972,0.7937); rgb(13pt)=(0.2661,0.2011,0.8043); rgb(14pt)=(0.2676,0.2052,0.8148); rgb(15pt)=(0.2691,0.2094,0.8249); rgb(16pt)=(0.2704,0.2138,0.8346); rgb(17pt)=(0.2717,0.2184,0.8439); rgb(18pt)=(0.2729,0.2231,0.8528); rgb(19pt)=(0.274,0.228,0.8612); rgb(20pt)=(0.2749,0.233,0.8692); rgb(21pt)=(0.2758,0.2382,0.8767); rgb(22pt)=(0.2766,0.2435,0.884); rgb(23pt)=(0.2774,0.2489,0.8908); rgb(24pt)=(0.2781,0.2543,0.8973); rgb(25pt)=(0.2788,0.2598,0.9035); rgb(26pt)=(0.2794,0.2653,0.9094); rgb(27pt)=(0.2798,0.2708,0.915); rgb(28pt)=(0.2802,0.2764,0.9204); rgb(29pt)=(0.2806,0.2819,0.9255); rgb(30pt)=(0.2809,0.2875,0.9305); rgb(31pt)=(0.2811,0.293,0.9352); rgb(32pt)=(0.2813,0.2985,0.9397); rgb(33pt)=(0.2814,0.304,0.9441); rgb(34pt)=(0.2814,0.3095,0.9483); rgb(35pt)=(0.2813,0.315,0.9524); rgb(36pt)=(0.2811,0.3204,0.9563); rgb(37pt)=(0.2809,0.3259,0.96); rgb(38pt)=(0.2807,0.3313,0.9636); rgb(39pt)=(0.2803,0.3367,0.967); rgb(40pt)=(0.2798,0.3421,0.9702); rgb(41pt)=(0.2791,0.3475,0.9733); rgb(42pt)=(0.2784,0.3529,0.9763); rgb(43pt)=(0.2776,0.3583,0.9791); rgb(44pt)=(0.2766,0.3638,0.9817); rgb(45pt)=(0.2754,0.3693,0.984); rgb(46pt)=(0.2741,0.3748,0.9862); rgb(47pt)=(0.2726,0.3804,0.9881); rgb(48pt)=(0.271,0.386,0.9898); rgb(49pt)=(0.2691,0.3916,0.9912); rgb(50pt)=(0.267,0.3973,0.9924); rgb(51pt)=(0.2647,0.403,0.9935); rgb(52pt)=(0.2621,0.4088,0.9946); rgb(53pt)=(0.2591,0.4145,0.9955); rgb(54pt)=(0.2556,0.4203,0.9965); rgb(55pt)=(0.2517,0.4261,0.9974); rgb(56pt)=(0.2473,0.4319,0.9983); rgb(57pt)=(0.2424,0.4378,0.9991); rgb(58pt)=(0.2369,0.4437,0.9996); rgb(59pt)=(0.2311,0.4497,0.9995); rgb(60pt)=(0.225,0.4559,0.9985); rgb(61pt)=(0.2189,0.462,0.9968); rgb(62pt)=(0.2128,0.4682,0.9948); rgb(63pt)=(0.2066,0.4743,0.9926); rgb(64pt)=(0.2006,0.4803,0.9906); rgb(65pt)=(0.195,0.4861,0.9887); rgb(66pt)=(0.1903,0.4919,0.9867); rgb(67pt)=(0.1869,0.4975,0.9844); rgb(68pt)=(0.1847,0.503,0.9819); rgb(69pt)=(0.1831,0.5084,0.9793); rgb(70pt)=(0.1818,0.5138,0.9766); rgb(71pt)=(0.1806,0.5191,0.9738); rgb(72pt)=(0.1795,0.5244,0.9709); rgb(73pt)=(0.1785,0.5296,0.9677); rgb(74pt)=(0.1778,0.5349,0.9641); rgb(75pt)=(0.1773,0.5401,0.9602); rgb(76pt)=(0.1768,0.5452,0.956); rgb(77pt)=(0.1764,0.5504,0.9516); rgb(78pt)=(0.1755,0.5554,0.9473); rgb(79pt)=(0.174,0.5605,0.9432); rgb(80pt)=(0.1716,0.5655,0.9393); rgb(81pt)=(0.1686,0.5705,0.9357); rgb(82pt)=(0.1649,0.5755,0.9323); rgb(83pt)=(0.161,0.5805,0.9289); rgb(84pt)=(0.1573,0.5854,0.9254); rgb(85pt)=(0.154,0.5902,0.9218); rgb(86pt)=(0.1513,0.595,0.9182); rgb(87pt)=(0.1492,0.5997,0.9147); rgb(88pt)=(0.1475,0.6043,0.9113); rgb(89pt)=(0.1461,0.6089,0.908); rgb(90pt)=(0.1446,0.6135,0.905); rgb(91pt)=(0.1429,0.618,0.9022); rgb(92pt)=(0.1408,0.6226,0.8998); rgb(93pt)=(0.1383,0.6272,0.8975); rgb(94pt)=(0.1354,0.6317,0.8953); rgb(95pt)=(0.1321,0.6363,0.8932); rgb(96pt)=(0.1288,0.6408,0.891); rgb(97pt)=(0.1253,0.6453,0.8887); rgb(98pt)=(0.1219,0.6497,0.8862); rgb(99pt)=(0.1185,0.6541,0.8834); rgb(100pt)=(0.1152,0.6584,0.8804); rgb(101pt)=(0.1119,0.6627,0.877); rgb(102pt)=(0.1085,0.6669,0.8734); rgb(103pt)=(0.1048,0.671,0.8695); rgb(104pt)=(0.1009,0.675,0.8653); rgb(105pt)=(0.0964,0.6789,0.8609); rgb(106pt)=(0.0914,0.6828,0.8562); rgb(107pt)=(0.0855,0.6865,0.8513); rgb(108pt)=(0.0789,0.6902,0.8462); rgb(109pt)=(0.0713,0.6938,0.8409); rgb(110pt)=(0.0628,0.6972,0.8355); rgb(111pt)=(0.0535,0.7006,0.8299); rgb(112pt)=(0.0433,0.7039,0.8242); rgb(113pt)=(0.0328,0.7071,0.8183); rgb(114pt)=(0.0234,0.7103,0.8124); rgb(115pt)=(0.0155,0.7133,0.8064); rgb(116pt)=(0.0091,0.7163,0.8003); rgb(117pt)=(0.0046,0.7192,0.7941); rgb(118pt)=(0.0019,0.722,0.7878); rgb(119pt)=(0.0009,0.7248,0.7815); rgb(120pt)=(0.0018,0.7275,0.7752); rgb(121pt)=(0.0046,0.7301,0.7688); rgb(122pt)=(0.0094,0.7327,0.7623); rgb(123pt)=(0.0162,0.7352,0.7558); rgb(124pt)=(0.0253,0.7376,0.7492); rgb(125pt)=(0.0369,0.74,0.7426); rgb(126pt)=(0.0504,0.7423,0.7359); rgb(127pt)=(0.0638,0.7446,0.7292); rgb(128pt)=(0.077,0.7468,0.7224); rgb(129pt)=(0.0899,0.7489,0.7156); rgb(130pt)=(0.1023,0.751,0.7088); rgb(131pt)=(0.1141,0.7531,0.7019); rgb(132pt)=(0.1252,0.7552,0.695); rgb(133pt)=(0.1354,0.7572,0.6881); rgb(134pt)=(0.1448,0.7593,0.6812); rgb(135pt)=(0.1532,0.7614,0.6741); rgb(136pt)=(0.1609,0.7635,0.6671); rgb(137pt)=(0.1678,0.7656,0.6599); rgb(138pt)=(0.1741,0.7678,0.6527); rgb(139pt)=(0.1799,0.7699,0.6454); rgb(140pt)=(0.1853,0.7721,0.6379); rgb(141pt)=(0.1905,0.7743,0.6303); rgb(142pt)=(0.1954,0.7765,0.6225); rgb(143pt)=(0.2003,0.7787,0.6146); rgb(144pt)=(0.2061,0.7808,0.6065); rgb(145pt)=(0.2118,0.7828,0.5983); rgb(146pt)=(0.2178,0.7849,0.5899); rgb(147pt)=(0.2244,0.7869,0.5813); rgb(148pt)=(0.2318,0.7887,0.5725); rgb(149pt)=(0.2401,0.7905,0.5636); rgb(150pt)=(0.2491,0.7922,0.5546); rgb(151pt)=(0.2589,0.7937,0.5454); rgb(152pt)=(0.2695,0.7951,0.536); rgb(153pt)=(0.2809,0.7964,0.5266); rgb(154pt)=(0.2929,0.7975,0.517); rgb(155pt)=(0.3052,0.7985,0.5074); rgb(156pt)=(0.3176,0.7994,0.4975); rgb(157pt)=(0.3301,0.8002,0.4876); rgb(158pt)=(0.3424,0.8009,0.4774); rgb(159pt)=(0.3548,0.8016,0.4669); rgb(160pt)=(0.3671,0.8021,0.4563); rgb(161pt)=(0.3795,0.8026,0.4454); rgb(162pt)=(0.3921,0.8029,0.4344); rgb(163pt)=(0.405,0.8031,0.4233); rgb(164pt)=(0.4184,0.803,0.4122); rgb(165pt)=(0.4322,0.8028,0.4013); rgb(166pt)=(0.4463,0.8024,0.3904); rgb(167pt)=(0.4608,0.8018,0.3797); rgb(168pt)=(0.4753,0.8011,0.3691); rgb(169pt)=(0.4899,0.8002,0.3586); rgb(170pt)=(0.5044,0.7993,0.348); rgb(171pt)=(0.5187,0.7982,0.3374); rgb(172pt)=(0.5329,0.797,0.3267); rgb(173pt)=(0.547,0.7957,0.3159); rgb(175pt)=(0.5748,0.7929,0.2941); rgb(176pt)=(0.5886,0.7913,0.2833); rgb(177pt)=(0.6024,0.7896,0.2726); rgb(178pt)=(0.6161,0.7878,0.2622); rgb(179pt)=(0.6297,0.7859,0.2521); rgb(180pt)=(0.6433,0.7839,0.2423); rgb(181pt)=(0.6567,0.7818,0.2329); rgb(182pt)=(0.6701,0.7796,0.2239); rgb(183pt)=(0.6833,0.7773,0.2155); rgb(184pt)=(0.6963,0.775,0.2075); rgb(185pt)=(0.7091,0.7727,0.1998); rgb(186pt)=(0.7218,0.7703,0.1924); rgb(187pt)=(0.7344,0.7679,0.1852); rgb(188pt)=(0.7468,0.7654,0.1782); rgb(189pt)=(0.759,0.7629,0.1717); rgb(190pt)=(0.771,0.7604,0.1658); rgb(191pt)=(0.7829,0.7579,0.1608); rgb(192pt)=(0.7945,0.7554,0.157); rgb(193pt)=(0.806,0.7529,0.1546); rgb(194pt)=(0.8172,0.7505,0.1535); rgb(195pt)=(0.8281,0.7481,0.1536); rgb(196pt)=(0.8389,0.7457,0.1546); rgb(197pt)=(0.8495,0.7435,0.1564); rgb(198pt)=(0.86,0.7413,0.1587); rgb(199pt)=(0.8703,0.7392,0.1615); rgb(200pt)=(0.8804,0.7372,0.165); rgb(201pt)=(0.8903,0.7353,0.1695); rgb(202pt)=(0.9,0.7336,0.1749); rgb(203pt)=(0.9093,0.7321,0.1815); rgb(204pt)=(0.9184,0.7308,0.189); rgb(205pt)=(0.9272,0.7298,0.1973); rgb(206pt)=(0.9357,0.729,0.2061); rgb(207pt)=(0.944,0.7285,0.2151); rgb(208pt)=(0.9523,0.7284,0.2237); rgb(209pt)=(0.9606,0.7285,0.2312); rgb(210pt)=(0.9689,0.7292,0.2373); rgb(211pt)=(0.977,0.7304,0.2418); rgb(212pt)=(0.9842,0.733,0.2446); rgb(213pt)=(0.99,0.7365,0.2429); rgb(214pt)=(0.9946,0.7407,0.2394); rgb(215pt)=(0.9966,0.7458,0.2351); rgb(216pt)=(0.9971,0.7513,0.2309); rgb(217pt)=(0.9972,0.7569,0.2267); rgb(218pt)=(0.9971,0.7626,0.2224); rgb(219pt)=(0.9969,0.7683,0.2181); rgb(220pt)=(0.9966,0.774,0.2138); rgb(221pt)=(0.9962,0.7798,0.2095); rgb(222pt)=(0.9957,0.7856,0.2053); rgb(223pt)=(0.9949,0.7915,0.2012); rgb(224pt)=(0.9938,0.7974,0.1974); rgb(225pt)=(0.9923,0.8034,0.1939); rgb(226pt)=(0.9906,0.8095,0.1906); rgb(227pt)=(0.9885,0.8156,0.1875); rgb(228pt)=(0.9861,0.8218,0.1846); rgb(229pt)=(0.9835,0.828,0.1817); rgb(230pt)=(0.9807,0.8342,0.1787); rgb(231pt)=(0.9778,0.8404,0.1757); rgb(232pt)=(0.9748,0.8467,0.1726); rgb(233pt)=(0.972,0.8529,0.1695); rgb(234pt)=(0.9694,0.8591,0.1665); rgb(235pt)=(0.9671,0.8654,0.1636); rgb(236pt)=(0.9651,0.8716,0.1608); rgb(237pt)=(0.9634,0.8778,0.1582); rgb(238pt)=(0.9619,0.884,0.1557); rgb(239pt)=(0.9608,0.8902,0.1532); rgb(240pt)=(0.9601,0.8963,0.1507); rgb(241pt)=(0.9596,0.9023,0.148); rgb(242pt)=(0.9595,0.9084,0.145); rgb(243pt)=(0.9597,0.9143,0.1418); rgb(244pt)=(0.9601,0.9203,0.1382); rgb(245pt)=(0.9608,0.9262,0.1344); rgb(246pt)=(0.9618,0.932,0.1304); rgb(247pt)=(0.9629,0.9379,0.1261); rgb(248pt)=(0.9642,0.9437,0.1216); rgb(249pt)=(0.9657,0.9494,0.1168); rgb(250pt)=(0.9674,0.9552,0.1116); rgb(251pt)=(0.9692,0.9609,0.1061); rgb(252pt)=(0.9711,0.9667,0.1001); rgb(253pt)=(0.973,0.9724,0.0938); rgb(254pt)=(0.9749,0.9782,0.0872); rgb(255pt)=(0.9769,0.9839,0.0805)},
colorbar
]
\addplot [forget plot] graphics [xmin=0.5, xmax=1024.5, ymin=0.5, ymax=1024.5] {fullQ-1.png};
\end{axis}
\end{tikzpicture}%
}
    \end{subfigure}\hfill
    \begin{subfigure}[b]{0.33\linewidth}
        \resizebox{\linewidth}{!}{
%
%
\begin{tikzpicture}

\begin{axis}[%
width=4.01in,
height=3.566in,
at={(0.673in,0.481in)},
scale only axis,
point meta min=0.7,
point meta max=1,
axis on top,
xmin=0.5,
xmax=1024.5,
xtick={0, 200, 400, 600, 800, 1000},
y dir=reverse,
ymin=0.5,
ymax=1024.5,
axis background/.style={fill=white},
title style={font=\bfseries},
colormap={mymap}{[1pt] rgb(0pt)=(0.2422,0.1504,0.6603); rgb(1pt)=(0.2444,0.1534,0.6728); rgb(2pt)=(0.2464,0.1569,0.6847); rgb(3pt)=(0.2484,0.1607,0.6961); rgb(4pt)=(0.2503,0.1648,0.7071); rgb(5pt)=(0.2522,0.1689,0.7179); rgb(6pt)=(0.254,0.1732,0.7286); rgb(7pt)=(0.2558,0.1773,0.7393); rgb(8pt)=(0.2576,0.1814,0.7501); rgb(9pt)=(0.2594,0.1854,0.761); rgb(11pt)=(0.2628,0.1932,0.7828); rgb(12pt)=(0.2645,0.1972,0.7937); rgb(13pt)=(0.2661,0.2011,0.8043); rgb(14pt)=(0.2676,0.2052,0.8148); rgb(15pt)=(0.2691,0.2094,0.8249); rgb(16pt)=(0.2704,0.2138,0.8346); rgb(17pt)=(0.2717,0.2184,0.8439); rgb(18pt)=(0.2729,0.2231,0.8528); rgb(19pt)=(0.274,0.228,0.8612); rgb(20pt)=(0.2749,0.233,0.8692); rgb(21pt)=(0.2758,0.2382,0.8767); rgb(22pt)=(0.2766,0.2435,0.884); rgb(23pt)=(0.2774,0.2489,0.8908); rgb(24pt)=(0.2781,0.2543,0.8973); rgb(25pt)=(0.2788,0.2598,0.9035); rgb(26pt)=(0.2794,0.2653,0.9094); rgb(27pt)=(0.2798,0.2708,0.915); rgb(28pt)=(0.2802,0.2764,0.9204); rgb(29pt)=(0.2806,0.2819,0.9255); rgb(30pt)=(0.2809,0.2875,0.9305); rgb(31pt)=(0.2811,0.293,0.9352); rgb(32pt)=(0.2813,0.2985,0.9397); rgb(33pt)=(0.2814,0.304,0.9441); rgb(34pt)=(0.2814,0.3095,0.9483); rgb(35pt)=(0.2813,0.315,0.9524); rgb(36pt)=(0.2811,0.3204,0.9563); rgb(37pt)=(0.2809,0.3259,0.96); rgb(38pt)=(0.2807,0.3313,0.9636); rgb(39pt)=(0.2803,0.3367,0.967); rgb(40pt)=(0.2798,0.3421,0.9702); rgb(41pt)=(0.2791,0.3475,0.9733); rgb(42pt)=(0.2784,0.3529,0.9763); rgb(43pt)=(0.2776,0.3583,0.9791); rgb(44pt)=(0.2766,0.3638,0.9817); rgb(45pt)=(0.2754,0.3693,0.984); rgb(46pt)=(0.2741,0.3748,0.9862); rgb(47pt)=(0.2726,0.3804,0.9881); rgb(48pt)=(0.271,0.386,0.9898); rgb(49pt)=(0.2691,0.3916,0.9912); rgb(50pt)=(0.267,0.3973,0.9924); rgb(51pt)=(0.2647,0.403,0.9935); rgb(52pt)=(0.2621,0.4088,0.9946); rgb(53pt)=(0.2591,0.4145,0.9955); rgb(54pt)=(0.2556,0.4203,0.9965); rgb(55pt)=(0.2517,0.4261,0.9974); rgb(56pt)=(0.2473,0.4319,0.9983); rgb(57pt)=(0.2424,0.4378,0.9991); rgb(58pt)=(0.2369,0.4437,0.9996); rgb(59pt)=(0.2311,0.4497,0.9995); rgb(60pt)=(0.225,0.4559,0.9985); rgb(61pt)=(0.2189,0.462,0.9968); rgb(62pt)=(0.2128,0.4682,0.9948); rgb(63pt)=(0.2066,0.4743,0.9926); rgb(64pt)=(0.2006,0.4803,0.9906); rgb(65pt)=(0.195,0.4861,0.9887); rgb(66pt)=(0.1903,0.4919,0.9867); rgb(67pt)=(0.1869,0.4975,0.9844); rgb(68pt)=(0.1847,0.503,0.9819); rgb(69pt)=(0.1831,0.5084,0.9793); rgb(70pt)=(0.1818,0.5138,0.9766); rgb(71pt)=(0.1806,0.5191,0.9738); rgb(72pt)=(0.1795,0.5244,0.9709); rgb(73pt)=(0.1785,0.5296,0.9677); rgb(74pt)=(0.1778,0.5349,0.9641); rgb(75pt)=(0.1773,0.5401,0.9602); rgb(76pt)=(0.1768,0.5452,0.956); rgb(77pt)=(0.1764,0.5504,0.9516); rgb(78pt)=(0.1755,0.5554,0.9473); rgb(79pt)=(0.174,0.5605,0.9432); rgb(80pt)=(0.1716,0.5655,0.9393); rgb(81pt)=(0.1686,0.5705,0.9357); rgb(82pt)=(0.1649,0.5755,0.9323); rgb(83pt)=(0.161,0.5805,0.9289); rgb(84pt)=(0.1573,0.5854,0.9254); rgb(85pt)=(0.154,0.5902,0.9218); rgb(86pt)=(0.1513,0.595,0.9182); rgb(87pt)=(0.1492,0.5997,0.9147); rgb(88pt)=(0.1475,0.6043,0.9113); rgb(89pt)=(0.1461,0.6089,0.908); rgb(90pt)=(0.1446,0.6135,0.905); rgb(91pt)=(0.1429,0.618,0.9022); rgb(92pt)=(0.1408,0.6226,0.8998); rgb(93pt)=(0.1383,0.6272,0.8975); rgb(94pt)=(0.1354,0.6317,0.8953); rgb(95pt)=(0.1321,0.6363,0.8932); rgb(96pt)=(0.1288,0.6408,0.891); rgb(97pt)=(0.1253,0.6453,0.8887); rgb(98pt)=(0.1219,0.6497,0.8862); rgb(99pt)=(0.1185,0.6541,0.8834); rgb(100pt)=(0.1152,0.6584,0.8804); rgb(101pt)=(0.1119,0.6627,0.877); rgb(102pt)=(0.1085,0.6669,0.8734); rgb(103pt)=(0.1048,0.671,0.8695); rgb(104pt)=(0.1009,0.675,0.8653); rgb(105pt)=(0.0964,0.6789,0.8609); rgb(106pt)=(0.0914,0.6828,0.8562); rgb(107pt)=(0.0855,0.6865,0.8513); rgb(108pt)=(0.0789,0.6902,0.8462); rgb(109pt)=(0.0713,0.6938,0.8409); rgb(110pt)=(0.0628,0.6972,0.8355); rgb(111pt)=(0.0535,0.7006,0.8299); rgb(112pt)=(0.0433,0.7039,0.8242); rgb(113pt)=(0.0328,0.7071,0.8183); rgb(114pt)=(0.0234,0.7103,0.8124); rgb(115pt)=(0.0155,0.7133,0.8064); rgb(116pt)=(0.0091,0.7163,0.8003); rgb(117pt)=(0.0046,0.7192,0.7941); rgb(118pt)=(0.0019,0.722,0.7878); rgb(119pt)=(0.0009,0.7248,0.7815); rgb(120pt)=(0.0018,0.7275,0.7752); rgb(121pt)=(0.0046,0.7301,0.7688); rgb(122pt)=(0.0094,0.7327,0.7623); rgb(123pt)=(0.0162,0.7352,0.7558); rgb(124pt)=(0.0253,0.7376,0.7492); rgb(125pt)=(0.0369,0.74,0.7426); rgb(126pt)=(0.0504,0.7423,0.7359); rgb(127pt)=(0.0638,0.7446,0.7292); rgb(128pt)=(0.077,0.7468,0.7224); rgb(129pt)=(0.0899,0.7489,0.7156); rgb(130pt)=(0.1023,0.751,0.7088); rgb(131pt)=(0.1141,0.7531,0.7019); rgb(132pt)=(0.1252,0.7552,0.695); rgb(133pt)=(0.1354,0.7572,0.6881); rgb(134pt)=(0.1448,0.7593,0.6812); rgb(135pt)=(0.1532,0.7614,0.6741); rgb(136pt)=(0.1609,0.7635,0.6671); rgb(137pt)=(0.1678,0.7656,0.6599); rgb(138pt)=(0.1741,0.7678,0.6527); rgb(139pt)=(0.1799,0.7699,0.6454); rgb(140pt)=(0.1853,0.7721,0.6379); rgb(141pt)=(0.1905,0.7743,0.6303); rgb(142pt)=(0.1954,0.7765,0.6225); rgb(143pt)=(0.2003,0.7787,0.6146); rgb(144pt)=(0.2061,0.7808,0.6065); rgb(145pt)=(0.2118,0.7828,0.5983); rgb(146pt)=(0.2178,0.7849,0.5899); rgb(147pt)=(0.2244,0.7869,0.5813); rgb(148pt)=(0.2318,0.7887,0.5725); rgb(149pt)=(0.2401,0.7905,0.5636); rgb(150pt)=(0.2491,0.7922,0.5546); rgb(151pt)=(0.2589,0.7937,0.5454); rgb(152pt)=(0.2695,0.7951,0.536); rgb(153pt)=(0.2809,0.7964,0.5266); rgb(154pt)=(0.2929,0.7975,0.517); rgb(155pt)=(0.3052,0.7985,0.5074); rgb(156pt)=(0.3176,0.7994,0.4975); rgb(157pt)=(0.3301,0.8002,0.4876); rgb(158pt)=(0.3424,0.8009,0.4774); rgb(159pt)=(0.3548,0.8016,0.4669); rgb(160pt)=(0.3671,0.8021,0.4563); rgb(161pt)=(0.3795,0.8026,0.4454); rgb(162pt)=(0.3921,0.8029,0.4344); rgb(163pt)=(0.405,0.8031,0.4233); rgb(164pt)=(0.4184,0.803,0.4122); rgb(165pt)=(0.4322,0.8028,0.4013); rgb(166pt)=(0.4463,0.8024,0.3904); rgb(167pt)=(0.4608,0.8018,0.3797); rgb(168pt)=(0.4753,0.8011,0.3691); rgb(169pt)=(0.4899,0.8002,0.3586); rgb(170pt)=(0.5044,0.7993,0.348); rgb(171pt)=(0.5187,0.7982,0.3374); rgb(172pt)=(0.5329,0.797,0.3267); rgb(173pt)=(0.547,0.7957,0.3159); rgb(175pt)=(0.5748,0.7929,0.2941); rgb(176pt)=(0.5886,0.7913,0.2833); rgb(177pt)=(0.6024,0.7896,0.2726); rgb(178pt)=(0.6161,0.7878,0.2622); rgb(179pt)=(0.6297,0.7859,0.2521); rgb(180pt)=(0.6433,0.7839,0.2423); rgb(181pt)=(0.6567,0.7818,0.2329); rgb(182pt)=(0.6701,0.7796,0.2239); rgb(183pt)=(0.6833,0.7773,0.2155); rgb(184pt)=(0.6963,0.775,0.2075); rgb(185pt)=(0.7091,0.7727,0.1998); rgb(186pt)=(0.7218,0.7703,0.1924); rgb(187pt)=(0.7344,0.7679,0.1852); rgb(188pt)=(0.7468,0.7654,0.1782); rgb(189pt)=(0.759,0.7629,0.1717); rgb(190pt)=(0.771,0.7604,0.1658); rgb(191pt)=(0.7829,0.7579,0.1608); rgb(192pt)=(0.7945,0.7554,0.157); rgb(193pt)=(0.806,0.7529,0.1546); rgb(194pt)=(0.8172,0.7505,0.1535); rgb(195pt)=(0.8281,0.7481,0.1536); rgb(196pt)=(0.8389,0.7457,0.1546); rgb(197pt)=(0.8495,0.7435,0.1564); rgb(198pt)=(0.86,0.7413,0.1587); rgb(199pt)=(0.8703,0.7392,0.1615); rgb(200pt)=(0.8804,0.7372,0.165); rgb(201pt)=(0.8903,0.7353,0.1695); rgb(202pt)=(0.9,0.7336,0.1749); rgb(203pt)=(0.9093,0.7321,0.1815); rgb(204pt)=(0.9184,0.7308,0.189); rgb(205pt)=(0.9272,0.7298,0.1973); rgb(206pt)=(0.9357,0.729,0.2061); rgb(207pt)=(0.944,0.7285,0.2151); rgb(208pt)=(0.9523,0.7284,0.2237); rgb(209pt)=(0.9606,0.7285,0.2312); rgb(210pt)=(0.9689,0.7292,0.2373); rgb(211pt)=(0.977,0.7304,0.2418); rgb(212pt)=(0.9842,0.733,0.2446); rgb(213pt)=(0.99,0.7365,0.2429); rgb(214pt)=(0.9946,0.7407,0.2394); rgb(215pt)=(0.9966,0.7458,0.2351); rgb(216pt)=(0.9971,0.7513,0.2309); rgb(217pt)=(0.9972,0.7569,0.2267); rgb(218pt)=(0.9971,0.7626,0.2224); rgb(219pt)=(0.9969,0.7683,0.2181); rgb(220pt)=(0.9966,0.774,0.2138); rgb(221pt)=(0.9962,0.7798,0.2095); rgb(222pt)=(0.9957,0.7856,0.2053); rgb(223pt)=(0.9949,0.7915,0.2012); rgb(224pt)=(0.9938,0.7974,0.1974); rgb(225pt)=(0.9923,0.8034,0.1939); rgb(226pt)=(0.9906,0.8095,0.1906); rgb(227pt)=(0.9885,0.8156,0.1875); rgb(228pt)=(0.9861,0.8218,0.1846); rgb(229pt)=(0.9835,0.828,0.1817); rgb(230pt)=(0.9807,0.8342,0.1787); rgb(231pt)=(0.9778,0.8404,0.1757); rgb(232pt)=(0.9748,0.8467,0.1726); rgb(233pt)=(0.972,0.8529,0.1695); rgb(234pt)=(0.9694,0.8591,0.1665); rgb(235pt)=(0.9671,0.8654,0.1636); rgb(236pt)=(0.9651,0.8716,0.1608); rgb(237pt)=(0.9634,0.8778,0.1582); rgb(238pt)=(0.9619,0.884,0.1557); rgb(239pt)=(0.9608,0.8902,0.1532); rgb(240pt)=(0.9601,0.8963,0.1507); rgb(241pt)=(0.9596,0.9023,0.148); rgb(242pt)=(0.9595,0.9084,0.145); rgb(243pt)=(0.9597,0.9143,0.1418); rgb(244pt)=(0.9601,0.9203,0.1382); rgb(245pt)=(0.9608,0.9262,0.1344); rgb(246pt)=(0.9618,0.932,0.1304); rgb(247pt)=(0.9629,0.9379,0.1261); rgb(248pt)=(0.9642,0.9437,0.1216); rgb(249pt)=(0.9657,0.9494,0.1168); rgb(250pt)=(0.9674,0.9552,0.1116); rgb(251pt)=(0.9692,0.9609,0.1061); rgb(252pt)=(0.9711,0.9667,0.1001); rgb(253pt)=(0.973,0.9724,0.0938); rgb(254pt)=(0.9749,0.9782,0.0872); rgb(255pt)=(0.9769,0.9839,0.0805)},
colorbar
]
\addplot [forget plot] graphics [xmin=0.5, xmax=1024.5, ymin=0.5, ymax=1024.5] {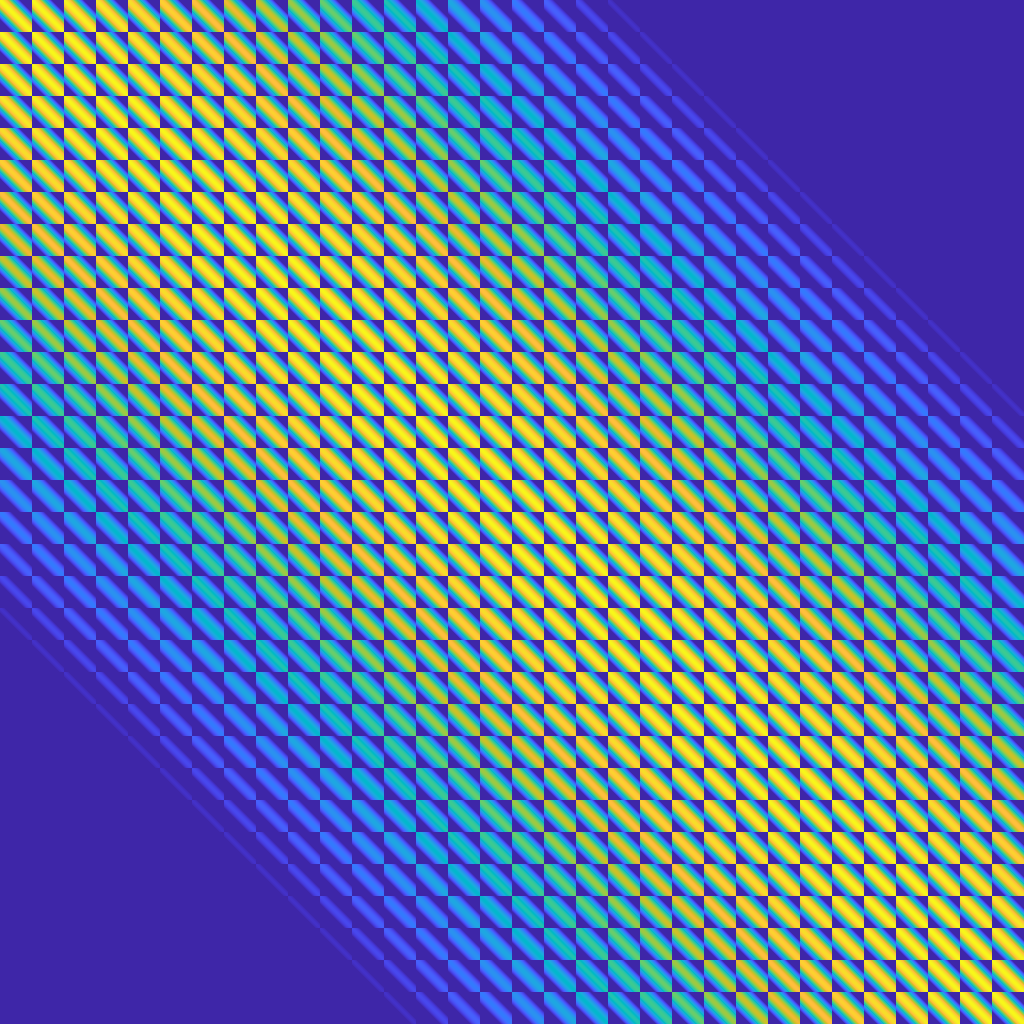};
\end{axis}
\end{tikzpicture}%
}
    \end{subfigure}\hfill
    \begin{subfigure}[b]{0.33\linewidth}
        \resizebox{\linewidth}{!}{\input{firstcol}}
    \end{subfigure}
    \caption{Left: the full Mat\'{e}rn covariance matrix $\mathbf{Q}$ ($\nu=1.5$, $\ell=0.9$) on a $32\times32$ grid. Center: the Kronecker--Toeplitz approximation $\hat{\mathbf{Q}}_1 = \mathbf{A}_r\otimes\mathbf{A}_c$ using the rank-1 PSF approximation. Right: the first columns of $\mathbf{Q}$ and $\hat{\mathbf{Q}}_1$ overlaid. We observe that the rank-1 approximation of the PSF captures the dominant structure of $\mathbf{Q}$.}
    \label{fig:origvsapprox}
\end{figure}

Consider 2D Mat\'{e}rn kernels $k_\nu(r)$ that lead to large covariance matrices $\mathbf{Q} \in \mathbb{R}^{n^2 \times n^2}$. 
To construct the Kronecker approximation, we take the first column of $\mathbf{Q}$, $\mathbf{q}_1 \in \mathbb{R}^{n^2}$, and reshape it into an $n \times n$ PSF matrix $\textbf{P} = \mat(\mathbf{q}_1)$~\cite{kamm2000}. See the PSF shown in the right plot of \cref{fig:kernelQpsf}. Throughout, $n$ is a power of two, and $\ell$ denotes the correlation length. Given the SVD of the PSF, $\textbf{P} = \mathbf{U}\bm{\Sigma} \mathbf{V}^T$ with $\mathbf{U}, \bm{\Sigma}, \mathbf{V} \in \mathbb{R}^{n \times n}$, we form
\begin{equation}
  \mathbf{c}_1 = \sqrt{\sigma_1}\, \mathbf{u}_1, \qquad \mathbf{r}_1 = \sqrt{\sigma_1}\, \mathbf{v}_1,
\end{equation}
where $\sigma_1$ is the largest singular value and $\mathbf{u}_1$ and $\mathbf{v}_1$ are the corresponding left and right singular vectors, i.e., the first columns of $\mathbf{U}$ and $\mathbf{V}$. The vectors $\mathbf{c}_1$ and $\mathbf{r}_1$ define the $n \times n$ Toeplitz matrices $\mathbf{A}_c$ and $\mathbf{A}_r$, respectively, which give the Kronecker--Toeplitz approximation,
\begin{equation}
  \hat{\mathbf{Q}}_1 = \mathbf{A}_r \otimes \mathbf{A}_c.
\end{equation}
Note that this approximation is built from the rank-1 approximation of the PSF, $\mathbf{c}_1 \mathbf{r}_1^\top$. More generally, using a rank-$k$ approximation of the PSF (e.g., corresponding to the first $k$ singular values of the PSF), we can get an approximation to $\mathbf{Q}$ given by
\begin{equation}
  \hat{\mathbf{Q}}_k = \sum_{i=1}^k \mathbf{A}_{r,i} \otimes \mathbf{A}_{c,i},
\end{equation}
where $\mathbf{A}_{c,i}$ and $\mathbf{A}_{r,i}$ are the Toeplitz matrices built from $\mathbf{c}_i = \sqrt{\sigma_i}\, \mathbf{u}_i$ and $\mathbf{r}_i = \sqrt{\sigma_i}\, \mathbf{v}_i$. Unless stated otherwise, $\hat{\mathbf{Q}}$ denotes $\hat{\mathbf{Q}}_1$. Each $\hat{\mathbf{Q}}_k$ is an $n^2 \times n^2$ matrix, and in all our experiments $\hat{\mathbf{Q}}_k$ is
SPD. In \Cref{fig:origvsapprox}, we provide matrices $\mathbf{Q}$ and $\hat{\mathbf{Q}}_1$ for a $32 \times 32$ grid, along with an overlay of plots of the first columns, which are nearly indistinguishable. Thus, we expect that approximations $\hat{\mathbf{Q}}_1$ capture the dominant structure and information from $\mathbf{Q}$, and hence can be used as good surrogates.

\paragraph{Storage and computational efficiency}

\begin{figure}[h!]
    \centering
        \includegraphics[width=1.0\linewidth]{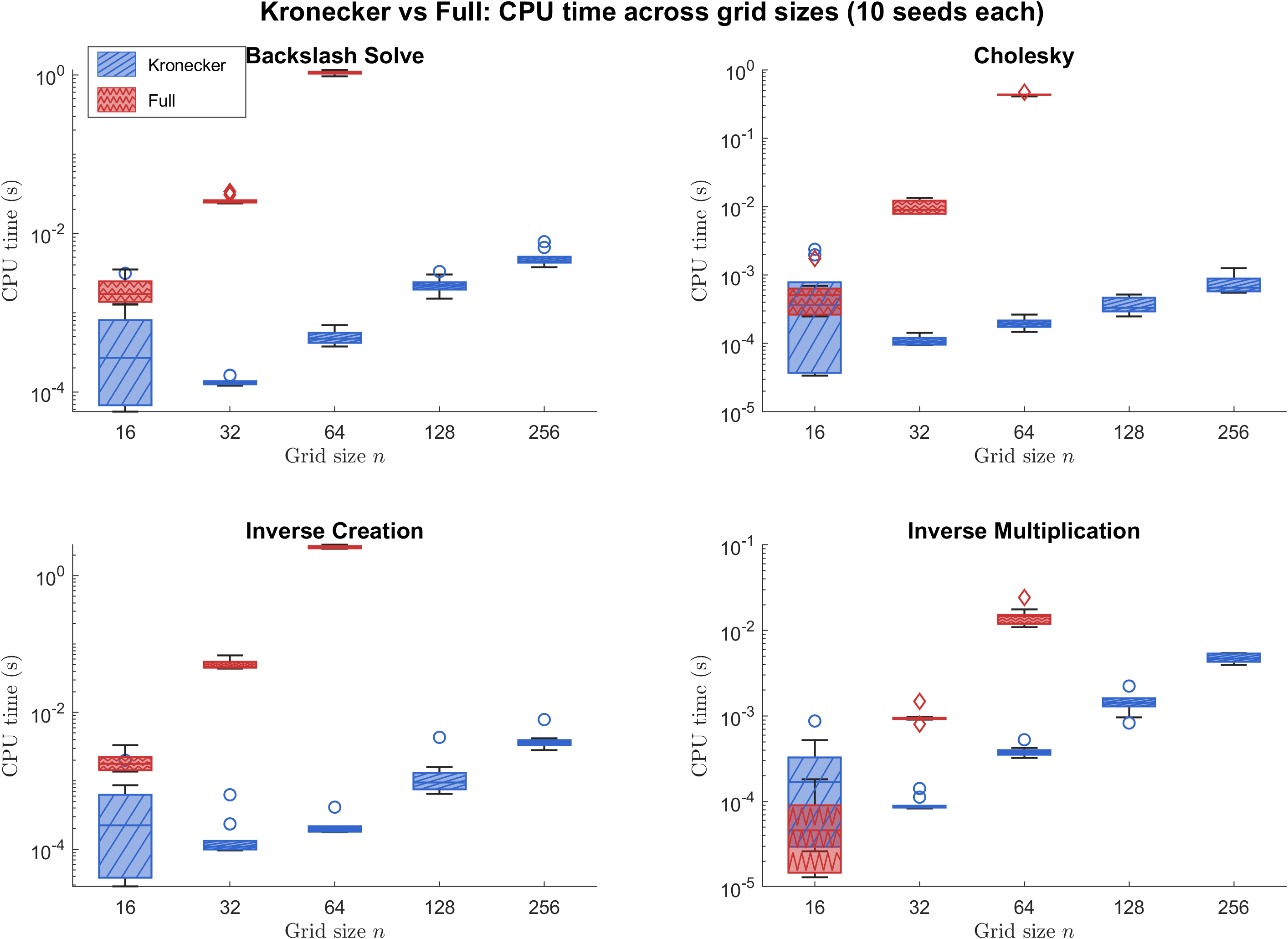}
        \caption{CPU times for operations with $\mathbf{Q}$ versus $\hat{\mathbf{Q}}_1$ for $\kappa(\nu = 1.5,\, \ell = 0.90)$ at several grid sizes over 10 random seeds. Each box spans the 25th to 75th percentiles with the median as the center line; whiskers extend to the most extreme runs within 1.5 times the interquartile range, and markers beyond them are outliers. The cost of operations involving $\mathbf{Q}$ increases rapidly with $n$, and $\mathbf{Q}$ becomes impractical for $n > 64$ as the memory requirement for its dense representation will be around 34.4GB at $n=256$. Operations with $\hat{\mathbf{Q}}_1$ depend on $n \times n$ Toeplitz factors and can be done in under $0.1$\,s.}
    \label{fig:diffseeds}
\end{figure}

For a $32 \times 32$ grid, the full matrix $\mathbf{Q} \in \mathbb{R}^{1024 \times 1024}$ requires 8.38\,MB of storage, whereas $\hat{\mathbf{Q}}$, stored as its two $32 \times 32$ Toeplitz factors, requires only 0.016\,MB, a reduction of roughly $500$ times. The Kronecker form also admits fast MVPs and inherits the algebraic properties of $\mathbf{Q}$; in particular, it is cheap to invert. Since $\hat{\mathbf{Q}}^{-1} = \mathbf{A}_r^{-1} \otimes \mathbf{A}_c^{-1}$, products with the inverse reuse the same Toeplitz factors, and by \eqref{eq:kronvec},
\begin{equation}
\label{eq:matfree}
  \hat{\mathbf{Q}}_1^{-1} \mathbf{x} = \vvec\!\left(\mathbf{A}_c^{-1} \mathbf{X} \mathbf{A}_r^{-T}\right), \qquad \mathbf{X} = \mat(\mathbf{x}) \in \mathbb{R}^{n \times n}.
\end{equation}
Applying $\hat{\mathbf{Q}}^{-1}$ thus requires only solves with two $n \times n$ Toeplitz matrices instead of one $n^2 \times n^2$ system. \Cref{fig:diffseeds} compares the CPU times of these operations with the corresponding operations with $\mathbf{Q}$ for different $n$. We observe that as the grid size grows, the CPU time for operations with the dense matrix $\mathbf{Q}$ increases rapidly, and $\mathbf{Q}$ becomes impractical for $n > 64$ due to its memory requirements. In contrast, operations with $\hat{\mathbf{Q}}_1$ involve only $n \times n$ Toeplitz factors and remain under $0.1$ s for all grid sizes tested, making the approximation practical for large-scale problems.

\subsection{Comparison with other matrix representations}
\label{sec:history}

We compare our Kronecker approximation with three alternative approaches in
\cref{sec:numcomp}. First, we consider a separable kernel representation, where the covariance matrix is a Kronecker product of two 1D Mat\`{e}rn covariance matrices.  Then we consider two approximations, including the RPCholesky method that gives a low-rank approximation of $\mathbf{Q}$ and the Poisson-based decomposition that gives an approximation of the \emph{inverse} $\mathbf{Q}^{-1}$. An approximation of $\mathbf{Q}$ can serve as a surrogate
(replacing $\mathbf{Q}$ in MVPs) or, if invertible, as a
preconditioner. 

\paragraph{Kronecker Product of Mat\'{e}rn Kernels} Covariance matrices with exact Kronecker product structure can arise naturally for dynamic inverse problems with spatiotemporal priors \cite{chung2018efficient} or for random fields with separable kernel functions. For example, Loh and Lam~\cite{loh2000} consider a Gaussian random field on $[0,1]^d$ with a separable (multiplicative) covariance, i.e., a product of one-dimensional kernels, one per coordinate direction. For $d = 2$ and unit variance, let $k(r)$ be a 1D kernel, e.g., \eqref{eq:matern}, and let $\{t_i\}_{i=1}^n$ be $n$ equispaced points in $[0,1]$. The 1D covariance matrix $\mathbf{Q}_{1D} \in \mathbb{R}^{n \times n}$ has elements
\begin{equation}
  [\mathbf{Q}_{1D}]_{ij} = k\!\left(|t_i - t_j|\right), \qquad i,j = 1, \ldots, n,
\end{equation}
and the separable 2D covariance matrix on the $n \times n$ grid is the Kronecker product
\begin{equation}
\label{eq:qkron}
  \mathbf{Q}_s = \mathbf{Q}_{1D} \otimes \mathbf{Q}_{1D} \in \mathbb{R}^{n^2 \times n^2},
\end{equation}
whose element for the grid nodes $(t_i, t_j)$ and $(t_k, t_l)$ is $k(|t_i - t_k|)\,k(|t_j - t_l|)$. The authors use this structure to obtain closed-form expressions for SPD matrix decompositions and inverses, with applications in the design and analysis of computer experiments. In general, the separable kernel is not the same as the isotropic 2D Mat\'{e}rn kernel \eqref{eq:matern}. The two coincide for the squared-exponential (Gaussian) kernel since $\exp(-(r_1^2 + r_2^2)/(2\ell^2)) = \exp(-r_1^2/(2\ell^2))\exp(-r_2^2/(2\ell^2))$.

\paragraph{Randomly Pivoted Cholesky (RPCholesky)} Chen et al.~\cite{chen2024} propose a fast randomized algorithm for low-rank approximation of a positive semidefinite matrix. It uses a pivoted partial Cholesky factorization in which each pivot
column is chosen by a randomized selection rule, producing the low-rank approximation $\hat{\mathbf{Q}}_{\mathrm{RPC}} = \mathbf{F} \mathbf{F}^T$ where $\mathbf{F} \in \mathbb{R}^{n^2 \times k}$ is obtained from $kn^2$ entry evaluations. The output is a rank-$k$ column Nystr\"om approximation, which was investigated for several kernel methods~\cite{halko2011, martinsson2020, troppwebber2023}. 
Since $\hat{\mathbf{Q}}_{\mathrm{RPC}}$ is rank deficient, using it directly as a preconditioner is not feasible; 
instead, it serves as a surrogate for MVPs with
$\mathbf{Q}$. Although RPCholesky does not exploit Toeplitz or Kronecker structure, it is a natural comparison for the cumulative approximation $\hat{\mathbf{Q}}_k$ because it
also builds up a low-rank approximation incrementally.
\paragraph{Poisson-based preconditioner}
Saibaba, Chung, and Petroske~\cite{saibaba2020} investigated preconditioners for
Mat\'{e}rn covariance matrices based on the link between Mat\'{e}rn kernels and powers of the finite-difference Laplacian. For $\nu = 3/2$ and $\nu = 5/2$, the preconditioner is $(-\Delta)^\gamma$ with $\gamma = 1$ and $\gamma = 2$, respectively, where $\Delta$ is the Laplacian operator discretized using finite differences. For $\nu = 1/2$, the Cholesky factorization of $(-\Delta)$ is used. This construction relies on the fact that the inverse of an integral operator with a Mat\'{e}rn kernel is a fractional differential operator. Therefore, the matrix
$\mathbf{G} = (-\Delta)^\gamma$ approximates the \emph{precision matrix}
$\mathbf{Q}^{-1}$ rather than $\mathbf{Q}$. Besides being SPD with $\mathbf{G}\mathbf{Q} \approx \mathbf{I}$, $\mathbf{G}$ is sparse, which allows sparse MVPs and triangular
solves.

\subsection{Numerical Comparison of Approximations for a Fixed \texorpdfstring{$\mathbf{Q}$}{Q}}
\label{sec:numcomp}

For each of the considered approximations of $\mathbf{Q}$, we provide relative errors using the Frobenius norm and the trace:

\begin{itemize}
    \item \textbf{Frobenius norm}:
    \begin{equation}
    d_{\text{Frob}}(\mathbf{Q}, \hat{\mathbf{Q}}_k) = \frac{\left\| \mathbf{Q} - \hat{\mathbf{Q}}_k \right\|_F}{\|\mathbf{Q}\|_F} = \frac{\sqrt{\operatorname{tr}[(\mathbf{Q} - \hat{\mathbf{Q}}_k)^2]}}{\sqrt{\operatorname{tr}[\mathbf{Q}^2]}}
    \end{equation}
    \item \textbf{Trace error}:
    \begin{equation}
    \label{eq:trace}
    d_{\text{Trace}}(\mathbf{Q}, \hat{\mathbf{Q}}_k) = \frac{\operatorname{tr}(\mathbf{Q}) - \operatorname{tr}(\hat{\mathbf{Q}}_k)}{\operatorname{tr}(\mathbf{Q})}
    \end{equation}
\end{itemize}

\begin{figure}[h!]
  \centering
    \includegraphics[width=0.9\linewidth]{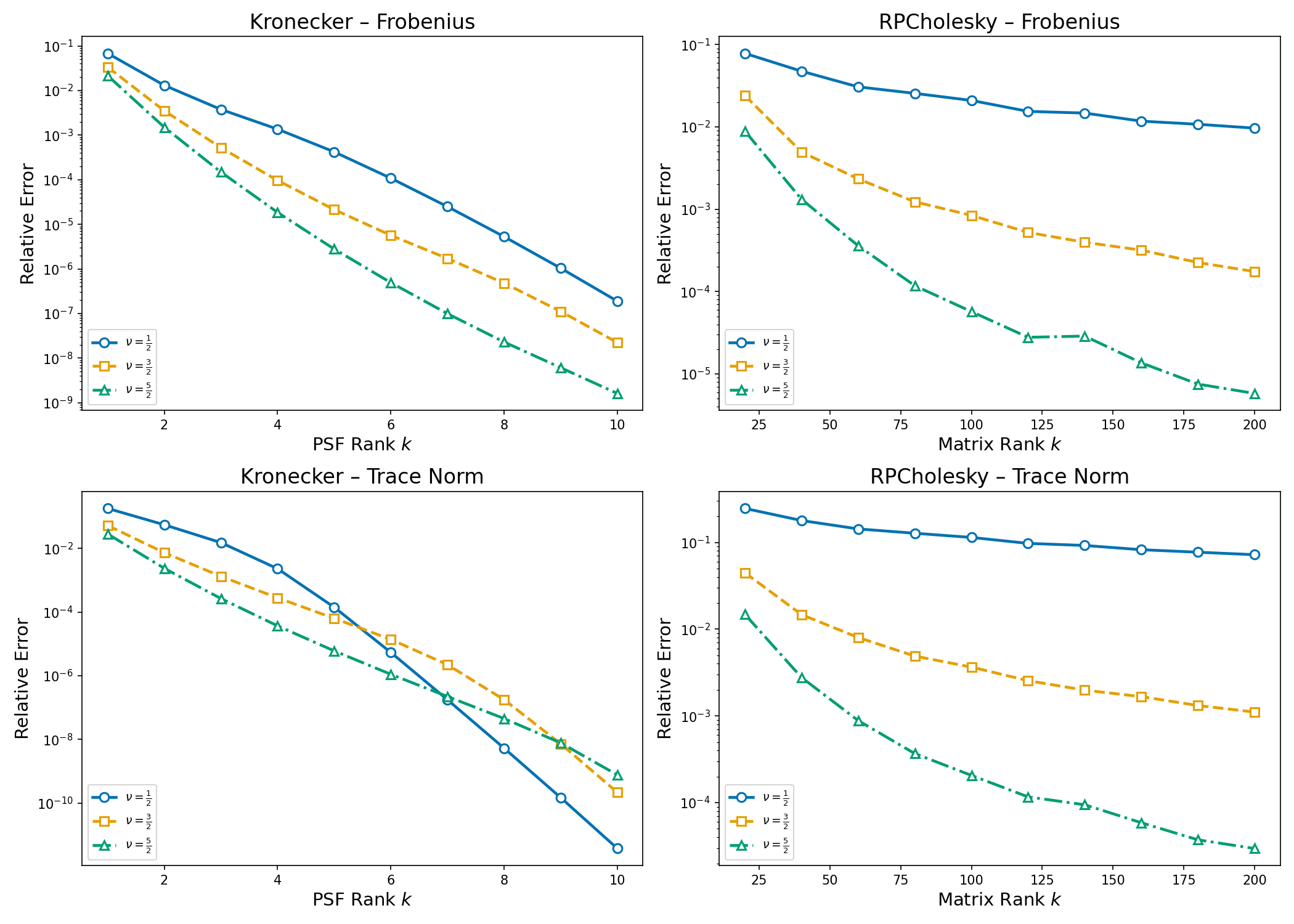}
    \caption{Relative Frobenius and trace errors of the Kronecker and RPCholesky approximations
    for Mat\'{e}rn kernels (left to right) with $\nu = \tfrac{1}{2}$, $\nu = \tfrac{3}{2}$, and $\nu = \tfrac{5}{2}$, all with $\ell = 0.50$.}
  \label{fig:matern-comparison}
\end{figure}

\begin{table}[h!]
\centering
\caption{Relative Frobenius and trace error norms between $\mathbf{Q}$ and its approximations for $\ell= 0.50$, $\nu = \frac12$ on a $32 \times 32$ grid, where the dense $\mathbf{Q}$ requires 8.38MB. 
All representations were converted to dense matrices before computing the errors. Storage for $\hat{\mathbf{Q}}$ counts the two $n\times n$ Toeplitz factors, and one $\textbf{Q}_{1D}$ factor for $\textbf{Q}_s$.}
\label{tab:method-compare}
\begin{tabular}{l|ccc}
\toprule
\textbf{Method} & \textbf{Frobenius} & \textbf{Trace Norm} & \textbf{Storage (MB)} \\
\midrule
$\hat{\mathbf{Q}}_1$ & 0.068 & 0.180 & 0.016 \\
$\hat{\mathbf{Q}}_2$ & 0.013 & 0.055 & 0.031 \\
$\hat{\mathbf{Q}}_3$ & 0.003 & 0.015 & 0.047 \\
$\hat{\mathbf{Q}}_{\text{RPC}}$ (rank-20) & 0.117 & 0.280 & 0.156 \\
${\mathbf{Q}}_s$ & 0.188 & 0 & 0.0078 \\
\bottomrule
\end{tabular}
\end{table}
Since $\mathbf{Q}$ and its separable Kronecker form $\mathbf{Q}_s$ both have unit diagonal, $d_{\text{Trace}}(\mathbf{Q},\mathbf{Q}_s)=0$ by \eqref{eq:trace} in \cref{tab:method-compare}. From the error decay curves in \Cref{fig:matern-comparison}, we observe that both the Kronecker and RPCholesky approximations show rapid error decay with PSF rank and matrix rank, respectively, but their relative performance depends on the smoothness parameter $\nu = \tfrac{1}{2}, \tfrac{3}{2}, \tfrac{5}{2}$ and on the norm. 

\begin{figure}[h!]
        \centering
        \resizebox{0.5\textwidth}{!}{\input{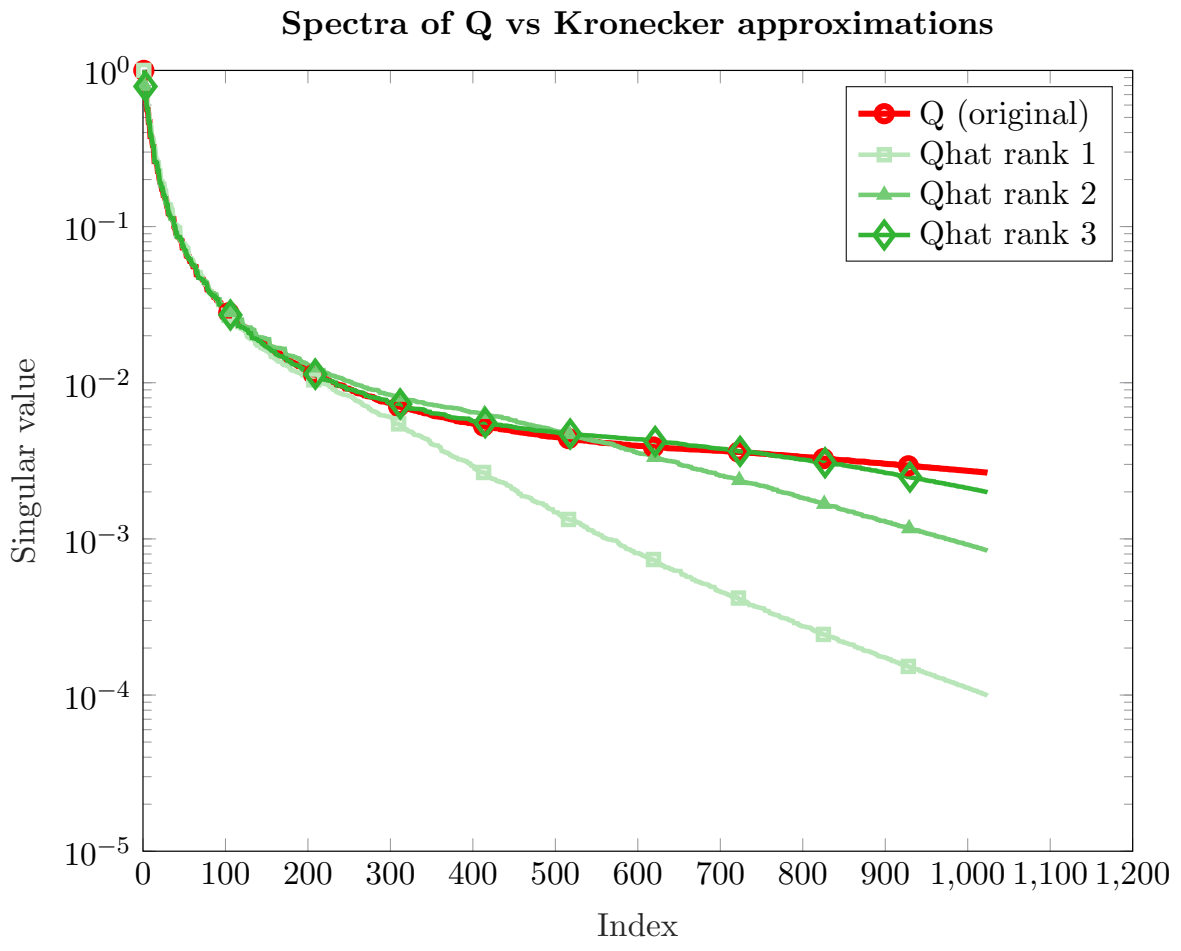}}
        \caption{Spectra of the SPD matrices $\hat{\mathbf{Q}}_k$ for various PSF ranks $k$. The Kronecker approximation captures the large singular values despite using only the dominant SVD components of the PSF. Beyond $\hat{\mathbf{Q}}_3$, higher PSF ranks $k$ give little improvement in spectral recovery.}
    \label{fig:fiveranks}
\end{figure}

Since $\mathbf{Q}$, its Kronecker--Toeplitz
approximation $\hat{\mathbf{Q}}$, and the Poisson-based preconditioner $\mathbf{G}$ are all SPD, we also measure the distances between them
with metrics specific to the geometry of SPD matrices. Let $\mathbf{Q}$ have the eigendecomposition $\mathbf{Q} = \mathbf{V}_Q \boldsymbol{\Lambda}_Q \mathbf{V}_Q^T$ with eigenvalues
$\lambda_\mathbf{Q} = (\lambda_1^\mathbf{Q}, \ldots, \lambda_{n^2}^\mathbf{Q})$ in non-increasing order, and define the eigenvalues $\lambda_k$ of $\hat{\mathbf{Q}}_k$ similarly. We consider the following metrics:
\begin{itemize}
\item \textbf{AIRM} (Affine-Invariant Riemannian Metric):
\begin{equation}
d_{\text{AIRM}}(\mathbf{Q}, \hat{\mathbf{Q}}_k) = \sqrt{\sum_{j=1}^{n^2} \left(\log \lambda_j^{\mathbf{Q}^{-1/2} \hat{\mathbf{Q}}_k \mathbf{Q}^{-1/2}}\right)^2}
\end{equation}
This is computed from the generalized eigenvalues of $\mathbf{L}_Q^{-1} \hat{\mathbf{Q}}_k \mathbf{L}_Q^{-T}$, where $\mathbf{L}_Q$ is the Cholesky factor of $\mathbf{Q}$.
\item \textbf{Jensen-Bregman LogDet}:
\begin{equation}
d_{\text{JBLogDet}}(\mathbf{Q}, \hat{\mathbf{Q}}_k) = \log \det\left(\frac{\mathbf{Q} + \hat{\mathbf{Q}}_k}{2}\right) - \frac{1}{2}\left(\log \det \mathbf{Q} + \log \det \hat{\mathbf{Q}}_k\right)
\end{equation}
\item \textbf{Spectral Log-Eigenvalue}:
\begin{equation}
d_{\text{SpecLog}}(\mathbf{Q}, \hat{\mathbf{Q}}_k) = \sum_{j=1}^{n^2} \left(\log \lambda_j^\mathbf{Q} - \log \lambda_j^{\mathbf{Q}_k}\right)^2
\end{equation}
with eigenvalues sorted in descending order. 

The Log-Euclidean metric, not shown, gives similar results. \Cref{fig:spdmetrics} reports these distances for different PSF ranks $k = 1, \ldots, 5$, where the errors corresponding to the Poisson approximation are provided in the horizontal line for reference. We notice that, regardless of the metric used, the gap between $\hat{\mathbf{Q}}_k$ and $\mathbf{Q}$ becomes smaller when the PSF rank $k$ grows. Specifically, when $k=1$, the gaps are roughly of the same order of magnitude as those from the Poisson approximation, but for $k\geq 2$, the Kronecker approximation is closer to $\mathbf{Q}$.
\end{itemize}
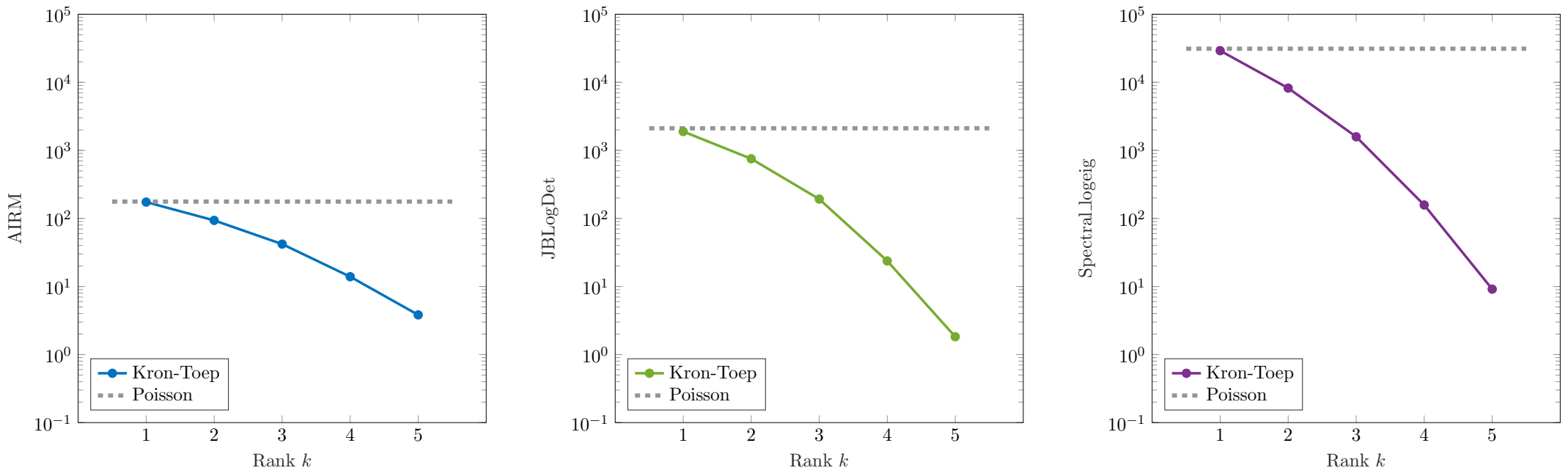
\begin{figure}[h!]
        \centering
        \resizebox{0.9\textwidth}{!}{
%
%
\definecolor{mycolor1}{rgb}{0.00000,0.44700,0.74100}%
\definecolor{mycolor2}{rgb}{0.46600,0.67400,0.18800}%
\definecolor{mycolor3}{rgb}{0.49400,0.18400,0.55600}%
\begin{tikzpicture}

\begin{axis}[%
width=3.334in,
height=3.334in,
at={(2.031in,1.298in)},
scale only axis,
xmin=0,
xmax=6,
xtick={1, 2, 3, 4, 5},
xlabel style={font=\color{white!15!black}},
xlabel={Rank $k$},
ymode=log,
ymin=0.1,
ymax=100000,
yminorticks=true,
ylabel style={font=\color{white!15!black}},
ylabel={AIRM},
axis background/.style={fill=white},
legend style={at={(0.03,0.03)}, anchor=south west, legend cell align=left, align=left, draw=white!15!black}
]
\addplot [color=mycolor1, line width=1.5pt, mark=*, mark options={solid, fill=mycolor1, mycolor1}]
  table[row sep=crcr]{%
1	173.933671122182\\
2	93.5485502221079\\
3	41.9660208097041\\
4	13.9311972218158\\
5	3.82599304080937\\
};
\addlegendentry{Kron-Toep}

\addplot [color=gray, dashed, line width=2.5pt]
  table[row sep=crcr]{%
0.5	176.734795220434\\
5.5	176.734795220434\\
};
\addlegendentry{Poisson}

\end{axis}

\begin{axis}[%
width=3.334in,
height=3.334in,
at={(6.419in,1.298in)},
scale only axis,
xmin=0,
xmax=6,
xtick={1, 2, 3, 4, 5},
xlabel style={font=\color{white!15!black}},
xlabel={Rank $k$},
ymode=log,
ymin=0.1,
ymax=100000,
yminorticks=true,
ylabel style={font=\color{white!15!black}},
ylabel={JBLogDet},
axis background/.style={fill=white},
legend style={at={(0.03,0.03)}, anchor=south west, legend cell align=left, align=left, draw=white!15!black}
]
\addplot [color=mycolor2, line width=1.5pt, mark=*, mark options={solid, fill=mycolor2, mycolor2}]
  table[row sep=crcr]{%
1	1891.04749718085\\
2	752.551986045365\\
3	192.36015220384\\
4	23.7349391085991\\
5	1.82647616268196\\
};
\addlegendentry{Kron-Toep}

\addplot [color=gray, dashed, line width=2.5pt]
  table[row sep=crcr]{%
0.5	2109.6556223519\\
5.5	2109.6556223519\\
};
\addlegendentry{Poisson}

\end{axis}

\begin{axis}[%
width=3.334in,
height=3.334in,
at={(10.806in,1.298in)},
scale only axis,
xmin=0,
xmax=6,
xtick={1, 2, 3, 4, 5},
xlabel style={font=\color{white!15!black}},
xlabel={Rank $k$},
ymode=log,
ymin=0.1,
ymax=100000,
yminorticks=true,
ylabel style={font=\color{white!15!black}},
ylabel={Spectral\_logeig},
axis background/.style={fill=white},
legend style={at={(0.03,0.03)}, anchor=south west, legend cell align=left, align=left, draw=white!15!black}
]
\addplot [color=mycolor3, line width=1.5pt, mark=*, mark options={solid, fill=mycolor3, mycolor3}]
  table[row sep=crcr]{%
1	29146.8832609915\\
2	8217.86436219706\\
3	1583.10457851442\\
4	157.550489421443\\
5	9.13623835930488\\
};
\addlegendentry{Kron-Toep}

\addplot [color=gray, dashed, line width=2.5pt]
  table[row sep=crcr]{%
0.5	31170.7954485526\\
5.5	31170.7954485526\\
};
\addlegendentry{Poisson}

\end{axis}
\end{tikzpicture}%
}
        \caption{Distance metrics for SPD matrices comparing $\hat{\mathbf{Q}}_k$ (rank-$k$ PSF Kronecker approximation, $k = 1, \ldots, 5$) and $\mathbf{Q}$ for $(\nu = 0.5, \ell = 0.50)$. Lower values indicate a better approximation. For reference, we provide distance metrics for the Poisson-based approximation, where $\textbf{G}^{-1}\approx \textbf{Q}$.}
    \label{fig:spdmetrics}
\end{figure}
Next we investigate the approximation of the spectra of $\mathbf{Q}_k$ for various PSF ranks $k$. As evident in \cref{fig:fiveranks}, $k=1$ is sufficient to capture the dominant eigenvalues of $\mathbf{Q},$ but a higher rank ($k=3$) is needed to approximate the entire spectrum of $\mathbf{Q}$ well. The required PSF rank depends on the purpose of the Kronecker approximation.  As we will show in \cref{sec:applications}, $k=1$ can be sufficient for solving an inverse problem, although higher PSF ranks can result in better solutions and can be easily incorporated in hybrid iterative methods for solving inverse problems. 

Across the aforementioned comparisons, the Kronecker approximations achieve lower relative Frobenius and trace norm errors, better spectral recovery, and greater SPD similarity to $\mathbf{Q}$. From here on, we consider only variants of the Kronecker approximation $\hat{\mathbf{Q}}$.

\section{Bayesian Inverse Problem Formulation and Iterative Solvers}
\label{sec:inverse} In this section, we describe how Kronecker-Toeplitz approximations described in \Cref{sec:approx} can be used to compute the MAP estimate efficiently for inverse problems. Consider again $\mathbf{b} = \mathbf{A}\mathbf{x} + \bm{\epsilon}$, where $\bm{\epsilon} \sim \mathcal{N}(\mathbf{0}, \mathbf{R})$ and $\mathbf{x} \sim \mathcal{N}(\mathbf{0}, \lambda^{-2}\mathbf{Q})$. By Bayes' rule, the posterior density can be written as $p(\mathbf{x} \mid \mathbf{b}) \propto p(\mathbf{b} \mid \mathbf{x}) \cdot p(\mathbf{x})$, with likelihood
\begin{equation}
  p(\mathbf{b} \mid \mathbf{x}) \propto \exp\!\left(-\frac{1}{2}(\mathbf{b} - \mathbf{A}\mathbf{x})^T \mathbf{R}^{-1}(\mathbf{b} - \mathbf{A}\mathbf{x})\right)
\end{equation}
and prior $p(\mathbf{x}) \propto \exp\!\left(-\frac{\lambda^2}{2}\, \mathbf{x}^T \mathbf{Q}^{-1} \mathbf{x}\right)$. Ignoring constants, we get the log posterior,
\begin{align}
  \log p(\mathbf{x} \mid \mathbf{b}) & = -\frac{1}{2}(\mathbf{b} - \mathbf{A}\mathbf{x})^T \mathbf{R}^{-1}(\mathbf{b} - \mathbf{A}\mathbf{x}) - \frac{\lambda^2}{2}\, \mathbf{x}^T \mathbf{Q}^{-1} \mathbf{x}\\
  & = -\frac{1}{2}\!\left[\mathbf{x}^T\!\left(\mathbf{A}^T \mathbf{R}^{-1} \mathbf{A} + \lambda^2 \mathbf{Q}^{-1}\right)\!\mathbf{x}
  - 2\mathbf{x}^T \mathbf{A}^T \mathbf{R}^{-1} \mathbf{b}\right].
\end{align}
Therefore, $\mathbf{x} \mid \mathbf{b} \sim \mathcal{N}(\bm{\mu}_{\mathrm{post}}, \mathbf{Q}_{\mathrm{post}})$, where
\begin{align}
    \label{eq:uq}
  \mathbf{Q}_{\mathrm{post}} &= \left(\mathbf{A}^T \mathbf{R}^{-1} \mathbf{A} + \lambda^2 \mathbf{Q}^{-1}\right)^{-1}, \\
  \bm{\mu}_{\mathrm{post}} &= \mathbf{Q}_{\mathrm{post}} \mathbf{A}^T \mathbf{R}^{-1} \mathbf{b}.
\end{align}
We minimize the negative log posterior to obtain the maximum a posteriori (MAP) estimate, which gives a solution to \cref{eq:tikhonov}. All of these quantities require solving large SPD linear systems with $\mathbf{Q}$ and $(\mathbf{A}^T \mathbf{R}^{-1} \mathbf{A} + \lambda^2 \mathbf{Q}^{-1})$. To reduce this cost, we use Krylov subspace methods~\cite{brown2025, sabate2025, chunggazzola2024}.

The matrix $\mathbf{A} \in \mathbb{R}^{m \times n^2}$ can be very large and is generally not square, which makes direct methods based on the SVD or QR decomposition expensive, costing $\mathcal{O}(m n^4)$ operations and requiring explicit access to $\mathbf{A}$. Krylov subspace methods avoid both issues~\cite{tropp2022, cortinovis2024, nakatsukasa2024}: using only MVPs with $\mathbf{A}$ and $\mathbf{A}^T$, they project the problem onto a small subspace built from those products. For rectangular $\mathbf{A}$, we use the Golub--Kahan (GK) bidiagonalization~\cite{GoVa13}, initialized with $\beta_1 \mathbf{u}_1 = \mathbf{b}$, where $\beta_1 = \|\mathbf{b}\|_2$, and $\alpha_1 \mathbf{v}_1 = \mathbf{A}^T \mathbf{u}_1$. After $j$ steps, GK produces two matrices with orthonormal columns,
\begin{equation*}
  \mathbf{U}_{j+1} = [\mathbf{u}_1, \ldots, \mathbf{u}_{j+1}] \in \mathbb{R}^{m \times (j+1)}, \qquad
  \mathbf{V}_j = [\mathbf{v}_1, \ldots, \mathbf{v}_j] \in \mathbb{R}^{n^2 \times j},
\end{equation*}
whose columns span the Krylov subspaces $\mathcal{K}_{j+1}(\mathbf{A}\mathbf{A}^T, \mathbf{b})$ and $\mathcal{K}_j(\mathbf{A}^T\mathbf{A}, \mathbf{A}^T\mathbf{b})$, respectively, where $\mathcal{K}_j(\mathbf{C}, \mathbf{g}) = \mathrm{span}\{\mathbf{g}, \mathbf{C}\mathbf{g}, \ldots, \mathbf{C}^{j-1}\mathbf{g}\}$. The matrices satisfy the GK relations
\begin{equation}
\label{eq:gk}
  \mathbf{A}\mathbf{V}_j = \mathbf{U}_{j+1}\mathbf{B}_j, \qquad
  \mathbf{A}^T\mathbf{U}_{j+1} = \mathbf{V}_j\mathbf{B}_j^T + \alpha_{j+1}\mathbf{v}_{j+1}\mathbf{e}_{j+1}^T,
\end{equation}
where $\mathbf{B}_j \in \mathbb{R}^{(j+1) \times j}$ is lower bidiagonal with diagonal entries $\alpha_1, \ldots, \alpha_j$ and subdiagonal entries $\beta_2, \ldots, \beta_{j+1}$. Setting $\mathbf{x} = \mathbf{V}_j \mathbf{y}$ with $\mathbf{y} \in \mathbb{R}^j$ and using $\mathbf{b} = \beta_1 \mathbf{U}_{j+1}\mathbf{e}_1$ reduces the least-squares problem to
\begin{equation}
  \|\mathbf{A}\mathbf{x} - \mathbf{b}\|_2 = \|\mathbf{U}_{j+1}(\mathbf{B}_j \mathbf{y} - \beta_1 \mathbf{e}_1)\|_2 = \|\mathbf{B}_j \mathbf{y} - \beta_1 \mathbf{e}_1\|_2,
\end{equation}
since $\mathbf{U}_{j+1}$ has orthonormal columns. The result is a $(j+1) \times j$ least-squares problem; solving it at each step gives the LSQR method~\cite{PaigeSaunders1982}. The generalized GK bidiagonalization used in GenHyBR extends the bidiagonalization process \eqref{eq:gk} to include weighted inner products defined by $\mathbf{R}^{-1}$ and $\mathbf{Q}$~\cite{chung2017}. This will be further described below.

For solving SPD systems, we use the \emph{conjugate gradient} (CG) method. CG is the canonical Krylov method for an SPD system $\mathbf{C}\mathbf{x} = \mathbf{d}$ with $\mathbf{C} \in \mathbb{R}^{n^2 \times n^2}$: it minimizes the $\mathbf{C}$-norm of the error over the growing Krylov subspace $\mathcal{K}_j(\mathbf{C}, \mathbf{d})$~\cite{GoVa13}. Each iteration requires one MVP with $\mathbf{C}$ and a few vector operations. Convergence slows as the condition number of $\mathbf{C}$ grows, which is why a preconditioner $\mathbf{M}$ that approximates $\mathbf{C}$ can be important.

With $\mathbf{M} \approx \mathbf{C}$, preconditioned CG (PCG) implicitly applies CG with the M inner product to the better-conditioned system $\mathbf{M}^{-1}\mathbf{C}\mathbf{x} = \mathbf{M}^{-1}\mathbf{d}$, at the cost of one application of $\mathbf{M}^{-1}$ per iteration ~\cite{doi:10.1137/1.9780898718003}. \Cref{alg:simple-pcg} gives the PCG implementation used in our
experiments; $\mathbf{M} = \mathbf{I}$ recovers standard CG. $\mathbf{x}^{(j)}$ is the $j$th iterate. We track the relative residual $\|\mathbf{r}^{(j)}\|/\|\mathbf{d}\|$ and, when the true solution $\mathbf{x}_\text{true}$ is available, the relative error $\|\mathbf{x}^{(j)} - \mathbf{x}_{\mathrm{true}}\|/\|\mathbf{x}_{\mathrm{true}}\|$.
\begin{algorithm}[h!]
\caption{Simple CG: solve $\mathbf{C}\mathbf{x} = \mathbf{d}$ with preconditioner $\mathbf{M}$}
\label{alg:simple-pcg}
\begin{algorithmic}
\STATE{\textbf{Result:} $[\mathbf{x}, \sim]$ = \texttt{simple\_pcg}($\mathbf{C}$, $\mathbf{d}$, \texttt{tol}, \texttt{maxiter}, $\mathbf{M}$, $\mathbf{x}_{\mathrm{true}}$})
\STATE{Initialize $\mathbf{x}^{(0)} := 0$, \; $\mathbf{r}^{(0)} := \mathbf{d} - \mathbf{C}\mathbf{x}^{(0)}$, \; $\mathbf{z}^{(0)} := \mathbf{M}^{-1}\mathbf{r}^{(0)}$, \; $\mathbf{p}^{(0)} := \mathbf{z}^{(0)}$, \; $\rho_0 := \mathbf{r}^{(0)T} \mathbf{z}^{(0)}$}
\FOR{$j = 0, 1, \ldots, \texttt{maxiter}-1$}
  \STATE{$\mathbf{q}^{(j)} := \mathbf{C} \mathbf{p}^{(j)}$}
  \STATE{$\alpha_j := \rho_j / (\mathbf{p}^{(j)T} \mathbf{q}^{(j)})$}
  \STATE{$\mathbf{x}^{(j+1)} := \mathbf{x}^{(j)} + \alpha_j \mathbf{p}^{(j)}$}
  \STATE{$\mathbf{r}^{(j+1)} := \mathbf{r}^{(j)} - \alpha_j \mathbf{q}^{(j)}$}
  \STATE{record $\|\mathbf{r}^{(j+1)}\|_2/\|\mathbf{d}\|_2$}
  \IF{$\|\mathbf{r}^{(j+1)}\|_2/\|\mathbf{d}\|_2 < \texttt{tol}$} \STATE{\textbf{return} $\mathbf{x}^{(j+1)}$} \ENDIF
  \STATE{$\mathbf{z}^{(j+1)} := \mathbf{M}^{-1} \mathbf{r}^{(j+1)}$}
  \STATE{$\rho_{j+1} := \mathbf{r}^{(j+1)T} \mathbf{z}^{(j+1)}$; \quad $\beta_j := \rho_{j+1}/\rho_j$}
  \STATE{$\mathbf{p}^{(j+1)} := \mathbf{z}^{(j+1)} + \beta_j \mathbf{p}^{(j)}$}
\ENDFOR
\end{algorithmic}
\end{algorithm}
Here, $\mathbf{C}$ is the SPD matrix or function handle, $\mathbf{d}$ the right-hand side, \texttt{tol} the convergence tolerance (e.g., $10^{-8}$), \texttt{maxiter} the maximum number of iterations (e.g., 100), $\mathbf{M}$ the preconditioner (default $\mathbf{M} = \mathbf{I}$), and $\mathbf{x}_{\mathrm{true}}$ the true solution used to compute relative reconstruction errors. The vector $\mathbf{r}^{(j)}$ is the residual, and $\mathbf{z}^{(j)}$ is the preconditioned residual.


Generalized hybrid iterative methods (GenHyBR), introduced by Chung and
Saibaba~\cite{chung2017}, extend GK to a $\mathbf{Q}$-weighted form. Recalling \eqref{eq:uq}, the MAP estimate is given by
\begin{equation}
    \label{eq:map}
  \mathbf{x}^{\mathrm{MAP}} = \arg\min_\mathbf{x} \frac{1}{2}\|\mathbf{A}\mathbf{x} - \mathbf{b}\|_{\mathbf{R}^{-1}}^2
  + \frac{\lambda^2}{2}\|\mathbf{x}\|_{\mathbf{Q}^{-1}}^2
  = \left(\mathbf{A}^T \mathbf{R}^{-1} \mathbf{A} + \lambda^2 \mathbf{Q}^{-1}\right)^{-1} \mathbf{A}^T \mathbf{R}^{-1} \mathbf{b},
\end{equation}
with $\bm{\epsilon} \sim \mathcal{N}(0, \mathbf{R})$ and $\mathbf{x} \sim \mathcal{N}(0, \lambda^{-2}\mathbf{Q})$. Standard hybrid solvers for general-form Tikhonov problems require a factorization $\mathbf{Q}^{-1} = \mathbf{L}_\mathbf{Q}^T \mathbf{L}_\mathbf{Q}$, known as ``priorconditioning'', which is infeasible when $\mathbf{Q}$ is a large, dense Mat\'{e}rn covariance matrix accessible only through MVPs. Chung and Saibaba~\cite{chung2017} show that, after a change of variables, $\mathbf{s} = \mathbf{Q}^{-1}\mathbf{x}$, we can solve a reformulated problem,
\begin{equation}
    \label{eq:qs}
  \min_\mathbf{s} \frac{1}{2}\|\mathbf{A}\mathbf{Q} \mathbf{s} - \mathbf{b}\|_{\mathbf{R}^{-1}}^2 + \frac{\lambda^2}{2}\|\mathbf{s}\|_{\mathbf{Q}}^2
\end{equation}
using the GenHyBR method.  These methods do not require directly computing the square root or inverse of the prior covariance matrix, and hence suit large-scale Bayesian inverse problems whose covariance kernels, such as those of the Mat\'{e}rn family, are accessible only through MVPs or are defined on irregular grids. The
resulting gen-LSQR iterates lie in the Krylov subspace $\mathcal{K}_j(\mathbf{A}^T \mathbf{R}^{-1} \mathbf{A} \mathbf{Q},\, \mathbf{A}^T \mathbf{R}^{-1} \mathbf{b})$, and the
shift invariance of this subspace in $\lambda$ allows the
regularization parameter to be updated for the projected problem at each iteration without rebuilding the basis.

At each GenHyBR iteration, a Tikhonov problem is solved on a lower-dimensional subspace, and the solution is then projected back to the original domain. Since the projected problem is small, $\lambda$ can be chosen at each step (e.g., by standard regularization parameter selection methods \cite{gazzola2019}). This combination of variational and iterative regularization is what makes GenHyBR a ``hybrid'' method. Each iteration requires one MVP with $\mathbf{A}$, one with $\mathbf{A}^T$, and two with $\mathbf{Q}$ (or $\hat{\mathbf{Q}}$).  By using \eqref{eq:kronvec}, MVPs with our Kronecker approximation $\hat{\mathbf{Q}}$ can be done efficiently.

The theory in~\cite{chung2017} also relates these methods to
classical approaches: for fixed $\lambda$ and in exact arithmetic, gen-LSQR iterates coincide with CG iterates (in the $\mathbf{Q}$-inner product) on the transformed normal equations. The iterates therefore inherit the favorable convergence of priorconditioning and the faster decay of the generalized singular values relative to those of $\mathbf{A}$, which allows early termination before noise dominates. 
\section{Using the Kronecker-Toeplitz Approximation for Solving Inverse Problems}
\label{sec:formulations}
Once $\hat{\mathbf{Q}}$ is constructed, it can be used to solve the inverse problem in several ways, depending on the computational goal. Our primary approach is to compute $\mathbf{x}_\text{MAP}$ from \eqref{eq:map} with $\mathbf{Q}$ replaced by $\hat{\mathbf{Q}}$ in different formulations.

\paragraph{$\hat{\mathbf{Q}}$ as a prior surrogate in GenHyBR}
Formulation \eqref{eq:qs} requires only MVPs with $\mathbf{Q}$, and GenHyBR selects the regularization parameter $\lambda$ automatically. Therefore, we propose to simply substitute $\hat{\mathbf{Q}}$ for $\mathbf{Q}$ and solve
\begin{equation}
  \min_\mathbf{s} \frac{1}{2}\|\mathbf{A}\hat{\mathbf{Q}} \mathbf{s} - \mathbf{b}\|_{\mathbf{R}^{-1}}^2 + \frac{\lambda^2}{2}\|\mathbf{s}\|_{\hat{\mathbf{Q}}}^2.
\end{equation}
\paragraph{$\hat{\mathbf{Q}}$ via CG on normal equations}
\label{sec:pcg}
A second approach applies standard, unpreconditioned CG to
\begin{equation}
    \label{eq:cg}
  \left(\mathbf{A}^T \mathbf{R}^{-1} \mathbf{A} + \lambda^2 \hat{\mathbf{Q}}^{-1}\right) \mathbf{x} = \mathbf{A}^T \mathbf{R}^{-1} \mathbf{b},
\end{equation}
where each CG iteration (\cref{alg:simple-pcg}) requires one MVP with $\mathbf{C} = \mathbf{A}^T \mathbf{R}^{-1} \mathbf{A} + \lambda^2 \hat{\mathbf{Q}}^{-1} \in \mathbb{R}^{n^2 \times n^2}$ for a fixed regularization parameter $\lambda$.
Substituting $\mathbf{x} = \hat{\mathbf{Q}}\mathbf{y}$ and left-multiplying by $\hat{\mathbf{Q}}$ gives
\begin{equation}
    \label{eq:reformedcg}
  \left(\hat{\mathbf{Q}}\mathbf{A}^T \mathbf{R}^{-1} \mathbf{A}\hat{\mathbf{Q}} + \lambda^2 \hat{\mathbf{Q}}\right) \mathbf{y} = \hat{\mathbf{Q}}\mathbf{A}^T \mathbf{R}^{-1} \mathbf{b},
\end{equation}
where we used the symmetry of $\hat{\mathbf{Q}}$. This formulation avoids $\hat{\mathbf{Q}}^{-1}$, but requires one to fix the regularization parameter in advance.

\begin{remark}
We also tested $\hat{\mathbf{Q}}$ as a preconditioner rather than a surrogate, so that when MVPs with $\mathbf{Q}^{-1}$ were needed in \eqref{eq:cg}, we used an inner CG solve, preconditioned by $\hat{\mathbf{Q}}$. This required two nested CG loops: the outer loop applies unpreconditioned CG to \eqref{eq:cg}, and the inner loop uses preconditioned CG for solves with $\mathbf{Q}$. For grids larger than $64 \times 64$, however, this approach was substantially slower than \eqref{eq:reformedcg} so we did not consider it further.
\end{remark}

\begin{remark} 
If the Cholesky factorization of the approximation $\hat{\mathbf{Q}}$ is readily available (e.g., $\hat{\mathbf{Q}_1}$, see \Cref{fig:diffseeds}), then one could use priorconditioning with the approximation and solve an augmented Tikhonov problem.  However, this is mathematically equivalent to applying CG to the normal equations \eqref{eq:cg} and does not provide a practical means for regularization parameter selection.
\end{remark}
 
\section{Numerical Experiments}
\label{sec:applications}
\label{sec:seismic} We consider the cross-well tomography problem from~\cite{chung2017}. A set of $s=256$ sources emit seismic signals sequentially, and another set of $p=50$ receivers record the travel time of each signal. This results in $m = s \cdot p = 256 \cdot 50 = 12,800$ observed travel times. We assume straight ray paths, a discretized domain in which each cell has constant slowness, and one measurement per source--receiver pair. Each travel time is the sum of the cell slownesses, weighted by the length of the ray segment in each cell. The goal is to recover the slowness of the medium from observed travel times that have been contaminated by Gaussian noise. The noise level is 4\% for all examples, i.e., $\|\epsilon\|_2 = 0.04\|\textbf{A}\textbf{x}_\text{true}\|_2$. Here, $\mathbf{A} \in \mathbb{R}^{12,800 \times 65,536}$ is the forward operator, $\mathbf{x}_\text{true} \in \mathbb{R}^{65,536}$ is the slowness on the $256 \times 256$ grid to be reconstructed, and $\mathbf{b} \in \mathbb{R}^{12,800}$ contains the Gaussian-noise-contaminated travel times, giving a problem of the form \cref{eq:noisyinv}. We further take $\mathbf{R} = \mathbf{I}$ as mentioned in \Cref{sec:intro}.

\begin{figure}[h!]
        \centering
        \resizebox{0.5\textwidth}{!}{\input{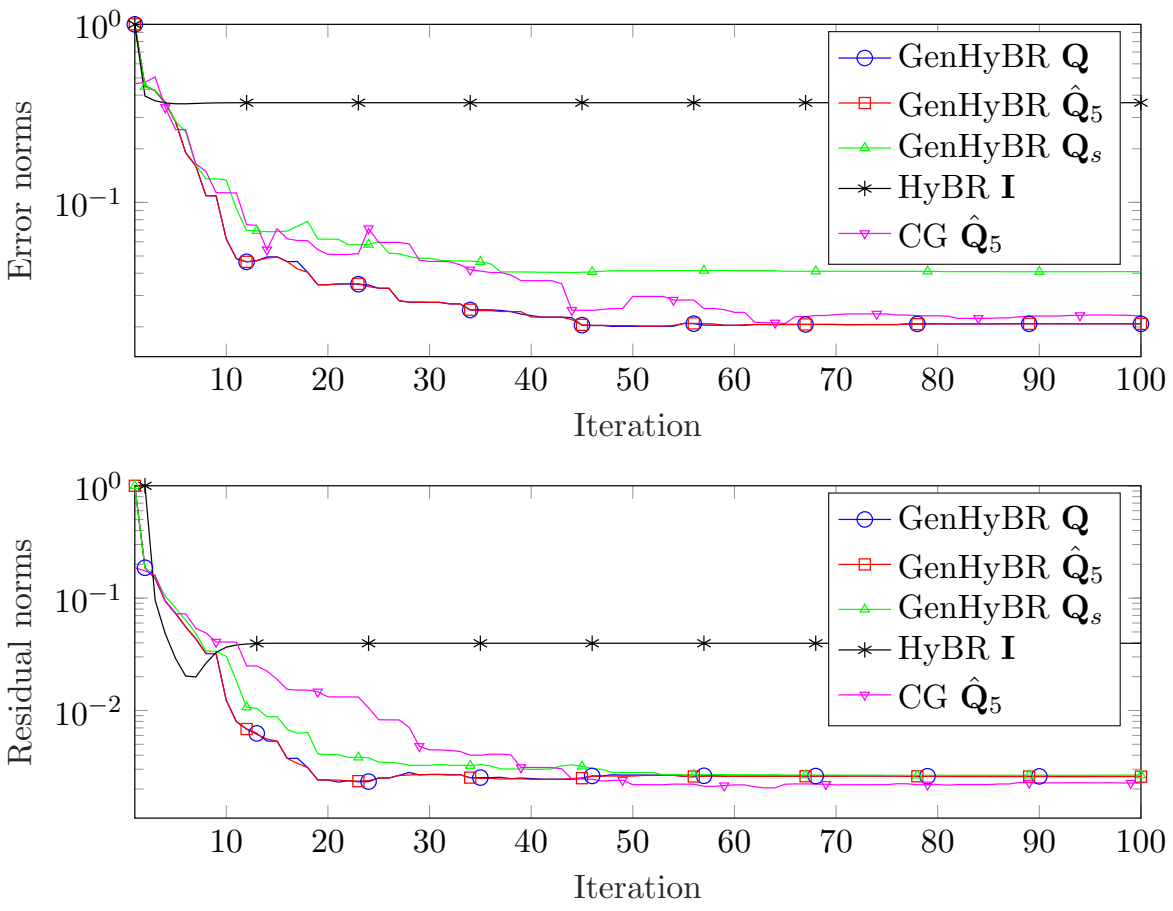}}
        \centering
    \caption{Relative reconstruction error $\|\mathbf{x}^{(j)} - \mathbf{x}_\mathrm{true}\| / \|\mathbf{x}_\mathrm{true}\|$ (top) and relative residual norm $\|\mathbf{A}\mathbf{x}^{(j)} - \mathbf{b}\| / \|\mathbf{b}\|$ (bottom) versus iteration for the seismic problem on a $256\times256$ grid with $\kappa(\nu=1.5, \ell=0.9)$. The method formulations include GenHyBR with the exact $\mathbf{Q}$, with the Kronecker approximation $\hat{\mathbf{Q}}$, and with the separable Kronecker form $\mathbf{Q}_s$, HyBR with $\mathbf{Q} = \mathbf{I}$, and CG on the reformulated normal equations. }
    \label{fig:threemethods}
\end{figure}

In this example, we seek an approximation of the MAP estimate \eqref{eq:map} where the prior covariance matrix $\mathbf{Q}$ represents a Matérn kernel $\kappa(\nu = 1.5, \ell=0.90)$. We compare various approaches and provide relative reconstruction errors and relative residual norms per iteration in \Cref{fig:threemethods}. We consider the GenHyBR method with $\mathbf{Q}$, $\hat{\mathbf{Q}}_5$, and $\mathbf{Q}_s$.  For comparison, we provide results obtained using HyBR, which corresponds to taking $\mathbf{Q} = \mathbf{I}$. In all of these methods, the optimal regularization parameter, i.e., the regularization parameter that minimizes the error between the reconstruction and the true solution, is used at each iteration. 
For a fixed regularization parameter, given by $\lambda^2 = 302.2787$, we also provide results for CG applied to the reformulated problem \eqref{eq:reformedcg}. The maximum number of iterations for all methods was set to 100.

We observe that error and residual curves for GenHyBR with $\hat{\mathbf{Q}}$ match those of GenHyBR with $\mathbf{Q}$, thus demonstrating that an approximate or surrogate prior covariance matrix can be used in the context of solving inverse problems.  
We observe that errors for HyBR are significantly larger than those of GenHyBR with $\mathbf{Q}$, revealing the importance of including an informative prior.  Moreover, errors for GenHyBR with $\mathbf{Q}_s$ level off near 0.05 instead of 0.03 for GenHyBR with $\hat{\mathbf{Q}}$, reflecting the difference in reconstructions that can be obtained with the choice of different priors.

    \begin{figure}[h!]
    \centering
        \resizebox{0.5\textwidth}{!}{
%
%
\begin{tikzpicture}

\begin{axis}[%
width=1.245in,
height=1.245in,
at={(0.758in,2.678in)},
scale only axis,
point meta min=0.00665123325459062,
point meta max=1,
axis on top,
xmin=0.5,
xmax=256.5,
tick align=outside,
y dir=reverse,
ymin=0.5,
ymax=256.5,
axis line style={draw=none},
ticks=none,
title style={font=\bfseries},
title={True}
]
\addplot [forget plot] graphics [xmin=0.5, xmax=256.5, ymin=0.5, ymax=256.5] {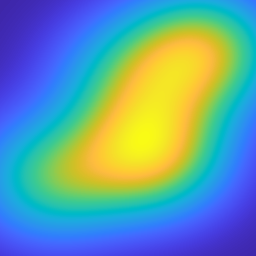};
\end{axis}

\begin{axis}[%
width=1.243in,
height=1.243in,
at={(2.396in,2.679in)},
scale only axis,
point meta min=-0.0336752358965344,
point meta max=0.99457990296939,
axis on top,
xmin=0.5,
xmax=256.5,
tick align=outside,
y dir=reverse,
ymin=0.5,
ymax=256.5,
axis line style={draw=none},
ticks=none,
title style={font=\bfseries},
title={GenHyBR Q}
]
\addplot [forget plot] graphics [xmin=0.5, xmax=256.5, ymin=0.5, ymax=256.5] {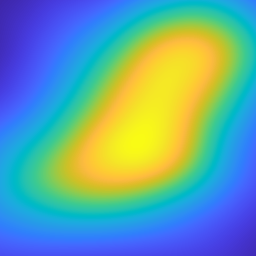};
\end{axis}

\begin{axis}[%
width=1.243in,
height=1.243in,
at={(4.034in,2.679in)},
scale only axis,
point meta min=-0.0371591365189311,
point meta max=0.976274418999489,
axis on top,
xmin=0.5,
xmax=256.5,
tick align=outside,
y dir=reverse,
ymin=0.5,
ymax=256.5,
axis line style={draw=none},
ticks=none,
title style={font=\bfseries},
title={GenHyBR $\hat{\mathbf{Q}}$}
]
\addplot [forget plot] graphics [xmin=0.5, xmax=256.5, ymin=0.5, ymax=256.5] {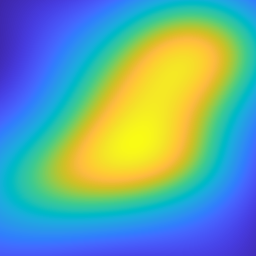};
\end{axis}

\begin{axis}[%
width=1.243in,
height=1.243in,
at={(0.758in,0.606in)},
scale only axis,
point meta min=-0.0606038542613054,
point meta max=0.980477567729789,
axis on top,
xmin=0.5,
xmax=256.5,
tick align=outside,
y dir=reverse,
ymin=0.5,
ymax=256.5,
axis line style={draw=none},
ticks=none,
title style={font=\bfseries},
title={GenHyBR Qs}
]
\addplot [forget plot] graphics [xmin=0.5, xmax=256.5, ymin=0.5, ymax=256.5] {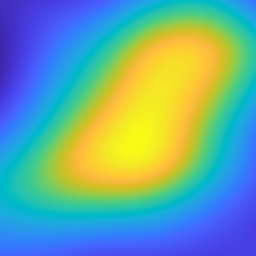};
\end{axis}

\begin{axis}[%
width=1.243in,
height=1.243in,
at={(2.396in,0.606in)},
scale only axis,
point meta min=-0.819817669633324,
point meta max=4.94708471845098,
axis on top,
xmin=0.5,
xmax=256.5,
tick align=outside,
y dir=reverse,
ymin=0.5,
ymax=256.5,
axis line style={draw=none},
ticks=none,
title style={font=\bfseries},
title={HyBR I}
]
\addplot [forget plot] graphics [xmin=0.5, xmax=256.5, ymin=0.5, ymax=256.5] {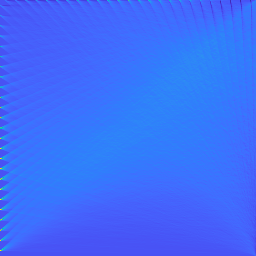};
\end{axis}

\begin{axis}[%
width=1.243in,
height=1.243in,
at={(4.034in,0.606in)},
scale only axis,
point meta min=-0.0419729847064195,
point meta max=1.01143892820758,
axis on top,
xmin=0.5,
xmax=256.5,
tick align=outside,
y dir=reverse,
ymin=0.5,
ymax=256.5,
axis line style={draw=none},
ticks=none,
title style={font=\bfseries},
title={CG $\hat{\mathbf{Q}}$}
]
\addplot [forget plot] graphics [xmin=0.5, xmax=256.5, ymin=0.5, ymax=256.5] {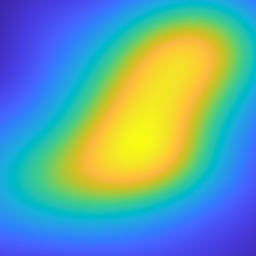};
\end{axis}
\end{tikzpicture}%
}
    \caption{True slowness field and MAP estimates on a $256\times256$ grid with $\kappa(\nu=1.5,\,\ell=0.90)$ for the same approaches as in \cref{fig:threemethods}: GenHyBR with $\mathbf{Q}$, $\hat{\mathbf{Q}}$, and ${\mathbf{Q}}_s$; HyBR with $\mathbf{I}$; and reformulated CG. Estimates based on the rank-5 Kronecker approximation $\hat{\mathbf{Q}}_5$ are visually indistinguishable from those with the exact $\mathbf{Q}$.}
    \label{fig:threemethodspic}
    \end{figure}
From \Cref{fig:threemethodspic}, we observe that the estimates obtained with GenHyBR with $\hat{\mathbf{Q}}_5$ and $\mathbf{Q}_s$, and with CG on the reformulated problem, closely match the estimate obtained with GenHyBR with $\mathbf{Q}$, whereas HyBR with $\mathbf{I}$ fails to recover the structure of the slowness field. In all cases, GenHyBR accesses $\hat{\mathbf{Q}}_k$ only through MVPs computed with \eqref{eq:kronvec}, so no $n^2 \times n^2$ matrix is ever formed.

\begin{figure}[h!]
    \centering
        \resizebox{0.9\textwidth}{!}{\input{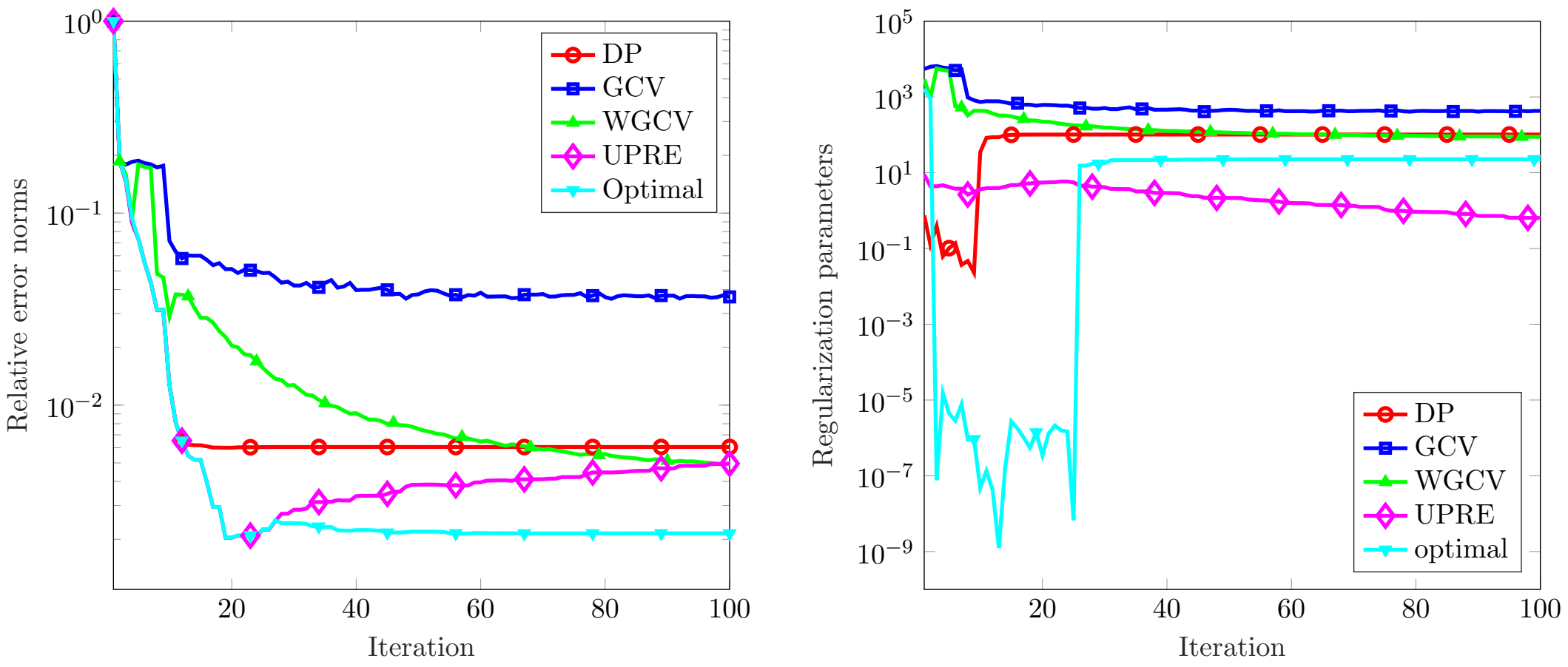}}
    \caption{Left: relative error norms for different regularization parameter choice methods in GenHyBR on a $32\times32$ grid with $\kappa(\nu=1.5, \ell = 0.90)$ and the Kronecker approximation $\hat{\mathbf{Q}}_5$; the original $\mathbf{Q}$ gives similar behavior for all methods. Right: computed regularization parameters per iteration. The methods are the discrepancy principle (DP), generalized cross-validation (GCV), weighted GCV (WGCV), the unbiased predictive risk estimator (UPRE), and the optimal regularization parameter. This corresponds to a noise standard deviation $\sigma = 0.6598$, which DP and UPRE use as the true $\sigma$.}
    \label{fig:threemethodspicregpar}
\end{figure}

\paragraph{Comparison of regularization parameter selection methods} GenHyBR selects a regularization parameter $\lambda_j$ at each iteration $j$ by applying a parameter choice rule to the small projected problem, so this selection is relatively inexpensive~\cite{chung2017, chung2024}. We consider five parameter choice rules all with $\hat{\mathbf{Q}}_5$. At each GenHyBR iteration $j$, $\lambda_j$ is chosen for the small projected problem, so all rules use only projected quantities: the projected residual $\mathbf{r}_j(\lambda)$, whose norm equals $\|\mathbf{b} - \mathbf{A}\mathbf{x}^{(j)}(\lambda)\|_2$, and the $(j+1) \times (j+1)$ projected influence matrix $\mathbf{H}_j(\lambda)$~\cite{chung2017}.
\begin{enumerate}
    \item The \emph{discrepancy principle} (DP) ~\cite{morozov1966} chooses $\lambda_j$ such that $\|\mathbf{r}_j(\lambda_j)\|_2^2 = \eta^2 m \sigma^2$, where $\sigma$ is the standard deviation of the noise and $\eta \geq 1$ is a safety factor; we use $\eta = 1$ and the true $\sigma$.
    \item The \emph{unbiased predictive risk estimator} (UPRE) ~\cite{turchin1967, vogel2002} minimizes $\|\mathbf{r}_j(\lambda)\|_2^2 + 2\sigma^2 \operatorname{trace}(\mathbf{H}_j(\lambda))$, an unbiased estimate of the predictive risk up to a constant, with the same $\sigma$.
    \item \emph{Generalized cross-validation} (GCV) ~\cite{golub1979} does not require the noise level and minimizes $\|\mathbf{r}_j(\lambda)\|_2^2 / \big(\operatorname{trace}(\mathbf{I}_{j+1} - \omega\mathbf{H}_j(\lambda))\big)^2$ with $\omega = 1$.
    \item \emph{Weighted GCV} (WGCV) ~\cite{chung2008} uses the same functional with a weight $\omega$ that is updated adaptively at each iteration to counteract the oversmoothing of GCV within hybrid methods~\cite{chung2024}.
    \item Finally, the optimal parameter minimizes $\|\mathbf{x}^{(j)}(\lambda) - \mathbf{x}_{\mathrm{true}}\|_2$; it requires $\mathbf{x}_{\mathrm{true}}$ and serves only as a reference.
\end{enumerate}
In \cref{fig:threemethodspicregpar}, for each of the considered regularization parameter choice methods, we provide the relative error norms per iteration and the computed values of $\lambda_j$ per iteration. We observe that relative error norms for DP level off after 12 iterations, thus demonstrating that an accurate noise estimate yields a good parameter early. GCV levels off near $4 \times 10^{-2}$, an order of magnitude above the other methods, because it selects parameters more than an order of magnitude larger than the optimal ones and therefore oversmooths. WGCV corrects this behavior, and its curve decreases steadily to about $5 \times 10^{-3}$ by iteration 100. UPRE follows the optimal curve until about iteration 25 and then drifts upward to about $5 \times 10^{-3}$, as its parameters fall below the optimal values.

\begin{figure}[h!]
        \centering
        \resizebox{0.5\textwidth}{!}{\input{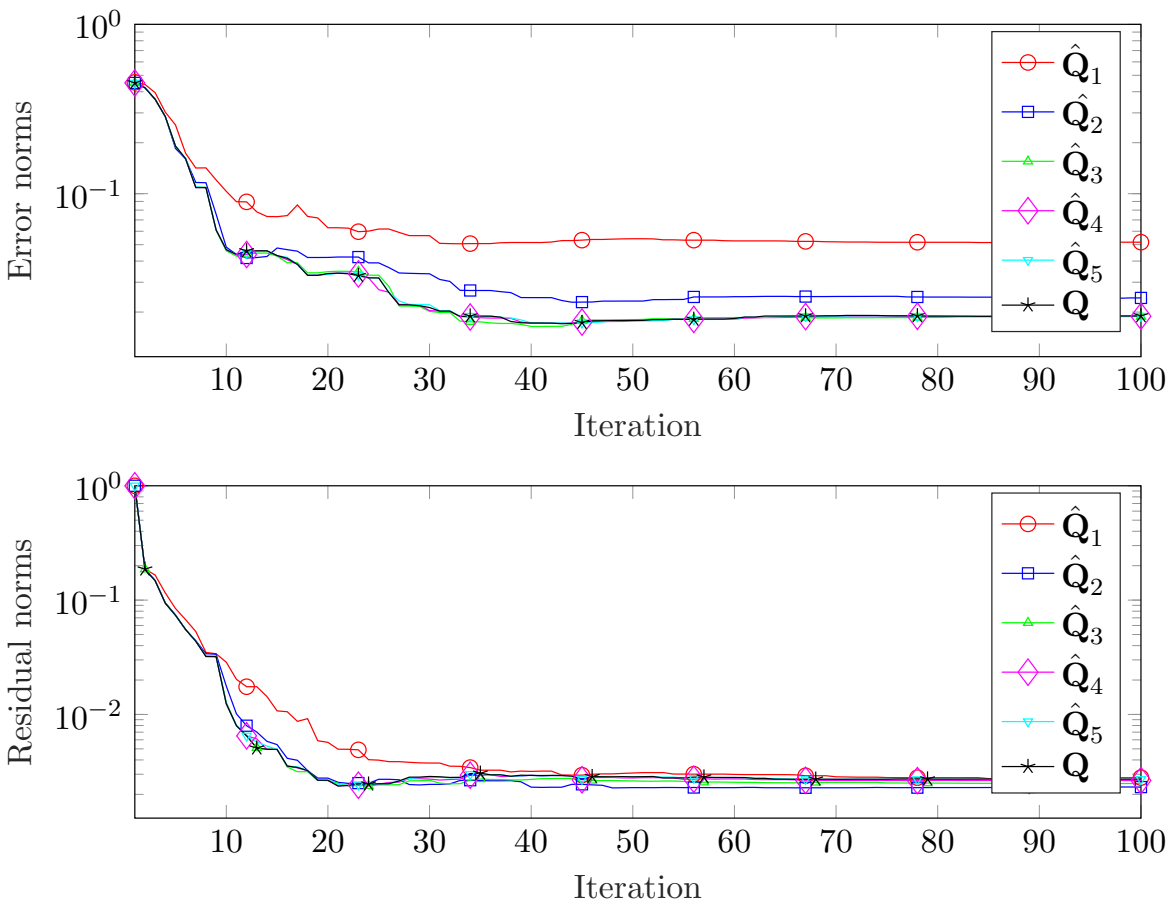}}
    \caption{Relative error and residual norms for GenHyBR with the cumulative Kronecker approximation $\hat{\mathbf{Q}}_k$, $k=1,\dots,5$, on a $256\times256$ grid with $\kappa(\nu=1.5,\,\ell=0.90)$ and the optimal regularization parameter. The curves for $k\geq3$ are virtually identical, indicating that $\hat{\textbf{Q}}_3$ is sufficient within GenHyBR, which accesses $\hat{\mathbf{Q}}_k$ only through MVPs.}
    \label{fig:GenHyBRranks}
\end{figure}

\paragraph{Investigating the impact of increasing PSF ranks in GenHyBR reconstructions} We repeat the experiment above using GenHyBR with the optimal regularization parameter and the cumulative approximation $\hat{\mathbf{Q}}_k = \sum_{i=1}^k \mathbf{A}_{r,i} \otimes \mathbf{A}_{c,i}$ for $k = 1, \ldots, 5$, and compare the results with GenHyBR with $\mathbf{Q}$ in \cref{fig:GenHyBRranks}. We observe that errors for $\hat{\mathbf{Q}}_1$ level off near 0.05, similar to those for $\mathbf{Q}_s$ in \cref{fig:threemethods}, whereas $\hat{\mathbf{Q}}_2$ already reduces the error to near 0.03, and the curves for $k \geq 3$ coincide with those for $\mathbf{Q}$. The residual norms show the same trend: $\hat{\mathbf{Q}}_1$ converges more slowly, but all approximations reach the same final residual level. Thus, a few Kronecker terms suffice to reproduce the reconstruction obtained with the exact $\mathbf{Q}$, and higher PSF ranks do not noticeably improve it. Since each additional term costs only two products of $n \times n$ matrices per MVP by \eqref{eq:kronvec}, increasing the rank adds little cost.

\section{Conclusions and Future Directions}
\label{sec:conclusions}
We investigated a Kronecker--Toeplitz approximation of large Mat\'{e}rn covariance matrices for Bayesian inverse problems. The approximation can replace the original covariance matrix directly in algorithms and iterative solvers. We validated both uses on a seismic tomography problem, using the GenHyBR framework to incorporate the Kronecker approximation into iterative methods. Compared with a low-rank approximation, a Poisson-based preconditioner, and a Kronecker product of two 1D Mat`{e}rn covariance matrices, the Kronecker-Toeplitz approximation achieved the lowest relative errors and closely captured the spectrum of the original covariance matrix. In future work, the approximation could be used in the posterior covariance \eqref{eq:uq} to enable efficient sampling without the full covariance matrix, and thus support uncertainty quantification. Upon publication, code will be available at \url{https://github.com/sarkarsrijon/kronapprox}.
\section*{Acknowledgments}
The author thanks Dr.\ Julianne Chung for guidance and support throughout this
research and acknowledges support from the Wilson Family Research Award.
This work was partially supported by the National Science Foundation under grant DMS-2411197 (J. Chung). Any opinions, findings, conclusions or recommendations expressed in this material are those of the author(s) and do not necessarily reflect the views of the National Science Foundation. We used large language models (Claude, Anthropic) to edit the manuscript for clarity and to help clean and test some numerical experiments. The authors reviewed and verified all text, code, and results, and take full responsibility for the content.

\bibliographystyle{siamplain}
\bibliography{references}

\end{document}